\documentclass[12pt,reqno]{amsart}
\usepackage{amssymb,amsmath,mathtools,amsthm,esint}
\usepackage{hyperref}
\usepackage[active]{srcltx}
\usepackage{fullpage}
\usepackage[T1]{fontenc}
\usepackage[ngerman]{babel}
\usepackage{xpatch}
\usepackage{algorithm}
\usepackage[n,advantage,operators,sets,adversary,landau,
probability,notions,logic,ff,mm,primitives,events,complexity,asymptotics,keys]{cryptocode}
\usepackage[noend]{algpseudocode}
\makeatletter
\def\BState{\State\hskip-\ALG@thistlm}
\makeatother
\xpatchcmd{\proof}{.}{\proofpunctuation}{}{}
\newcommand{\proofpunctuation}{:}
\title[Asymptotische Absch\"atzungen]
{Verteilung der Primzahlen und Asymptotik bez\"uglich der Kongruenzen der rationalen Punkte auf elliptischen Kurven \"uber einem endlichen K\"orper}
\author[Thi Altenschmidt]{Thi Altenschmidt}
\address{Fakult\"at f\"ur Elektrotechnik und Informationstechnik - IT-Sicherheit/Netze und Systeme - Ruhr Universit\"at Bochum}
\email{le.altenschmidt@rub.de}
\thanks{Ich m\"ochte mich an dieser Stelle bei meiner Ehefrau f\"ur ihre unendliche Liebe, Geduld und Unterst\"utzung bedanken.}
\begin{document}
\theoremstyle{plain}
\newtheorem{thm}{Satz}[section]
\newtheorem{satz}[thm]{Satz}
\newtheorem{lemma}[thm]{Lemma}
\newtheorem{folgerung}[thm]{Folgerung}
\newtheorem{hilfssatz}[thm]{Hilfssatz}
\newtheorem{bspl}[thm]{Beispiel}
\theoremstyle{definition}
\newtheorem{bezeichnungen}[thm]{Bezeichnungen}
\newtheorem{bem}[thm]{Bemerkung}
\numberwithin{equation}{thm}
\renewcommand{\theequation}{\arabic{section}.\arabic{thm}.\arabic{equation}}
\begin{abstract}
[English version] This work has two main purposes. On the one side we investigate in this work a question of H. Esnault on congruence formula in a construction of H. Esnault and C. Xu for the number of rational points on the closed fiber of a singular model of the projective plane over a local field. From the viewpoint of asymptotic analysis, the question is quite familiar with a question of N. Koblitz, which in turn has some meaningful applications in cryptography. We don't try to solve those questions in this work, but rather concentrate on studying asymptotic behaviours with very elementary techniques of generating functions. On the other side we extend the discussion on generating functions to the global situation, which inherit hybrid-properties from the Riemann zeta function and the $L$-function of elliptic curves. At the end we will look at an example of modular forms, where we are able to prove an analytic result, which is similar to a result of P\'olya for the Riemann $\xi$-function.    
\\
-----
\\
Wir befassen uns in dieser Arbeit zuerst mit einer Frage von H. Esnault \"uber die Kongruenzformel in einer Konstruktion von H. Esnault und C. Xu f\"ur die Anzahl der rationalen Punkte auf der abgeschlossenen Faser eines nicht-regul\"aren Modells von der projektiven Ebene \"uber einem lokalen K\"orper. Aus der Sicht der asymptotischen Analysis ist diese Frage ziemlich verwandt mit einem von N. Koblitz gestellten Problem, das wiederum eine relevante Bedeutung f\"ur die Anwendungen in der Kryptographie mit elliptischen Kurven hat. Dabei versuchen wir nicht gezielt, eine exakte geschlossene Kongruenzformel ausfindig zu machen, sondern gezielt, uns auf die asymptotische Analyse zu konzentrieren. Die zentrale Technik in unserer Arbeit ist die erzeugende Funktion aus der analytischen Kombinatorik, die wir im zweiten Teil dieser Arbeit auf den K\"orper der rationalen Zahlen \"ubertragen. Die erzeugende Funktion hat eine Hybrid-Form von der Riemannschen Zetafunktion und der $L$-Funktion der elliptischen Kurven. Am Ende betrachten wir als ein Bespiel die Spitzenformen vom Gewicht $2$ und beweisen ein analytisches Ergebnis, das \"ahnlich wie eine klassische Arbeit von P\'olya ist. 
\end{abstract}
\maketitle
\tableofcontents
\section*{Einleitung}
\noindent Sei $X/K$ eine projektive Variet\"at \"uber einem lokalen K\"orper $K$ mit einem endlichen Restk\"orper $\mathbb{F}_q$, wobei $q=p^n$ f\"ur eine Primzahl $p$ ist. Unter einem Modell von $X$ versteht man einen flachen projektiven Morphismus $\mathcal{X} \rightarrow \mathrm{Spec}\,\mathcal{O}_K$, so dass $\mathcal{X}\otimes_{\mathcal{O}_K}K \cong X$. Wenn $X/K$ glatt, projektiv und absolut irreduzible ist, haben H. Esnault und C. Xu in \cite[Thm. 1.1]{EsnXu09} bewiesen, dass jedes Modell $\mathcal{X}/\mathcal{O}_K$ mindestens einen $\mathbb{F}_q$-rationalen Punkt besitzt. Ferner, sie haben mithilfe von der Artinschen Konstruktion eines nicht-regul\"aren Modells $\mathcal{X}$ von $\mathbb{P}^2_K$ (vgl. \cite[\S 3]{EsnXu09}) konstruiert, so dass es 
\[\#\mathcal{X}(\mathbb{F}_q) \not\equiv 1 \mod q\]
gilt. Ihre Konstruktion zeigt, dass die Kongruenz der $\mathbb{F}_q$-rationalen Punkte der abgeschlossenen Faser abh\"angig von der Annahme der Regularit\"at des Modells ist. Denn H. Esnault hat bereits zuvor in \cite[Thm 1.1]{Esn06} gezeigt, dass falls $K$ ein $p$-adischer K\"orper ist und die allgemeine Faser $X/K$ bedingt auf der Coniveau-Filtrierung der $\ell$-adischen Kohomologie die Gleichheit
\[N^1H^i_{et}(X\otimes_KK^{alg},\mathbb{Q}_{\ell}) = H^i_{et}(X\otimes_KK^{alg},\mathbb{Q}_{\ell})\]
erf\"ullt, dann erf\"ullt jedes regul\"are Modell $\mathcal{X}/\mathcal{O}_K$ die Kongruenz
\[\#\mathcal{X}(\mathbb{F}_q) \equiv 1 \mod q,\]
wobei $K^{alg}$ ein algebraischer Abschluss von $K$ ist. Wir werden hierbei nicht \"uber die Maschinerie der $\ell$-adischen Kohomologie sowie Coniveau-Filtrierung oder die algebraische Geometrie der algebraischen R\"aume von Artin diskutieren. F\"ur den Rest in dieser Arbeit setzen wir nur Grundkenntnisse mit elliptischen Kurven \"uber endlichen K\"orpern und ein wenig \"uber den K\"orper der rationalen Zahlen $\mathbb{Q}$ voraus, die man leicht in \cite[\S V, \S VIII]{Sil09} nachschlagen kann, wenn es sich um die Geometrie handelt. Das Beispiel von H. Esnault und C. Xu (vgl. \cite[Thm. 1.2]{EsnXu09}) st\"utzt sich bis auf die sehr komplizierte geometrische Konstruktion auf die Tatsache, dass es zu jeder elliptischen Kurve $E/\mathbb{F}_q$ eine endliche K\"orpererweiterung $\mathbb{F}_{q^n}/\mathbb{F}_q$ mit
\[\#E(\mathbb{F}_{q^n}) \not\equiv 1 \mod q^n\]
gibt. In der Konstruktion von H. Esnault und C. Xu lautet die Kongruenzformel
\[\#\mathcal{X}(\mathbb{F}_q) \equiv 2 - \#E(\mathbb{F}_q) \mod q.\]
Um eine geschlossene Kongruenzformel f\"ur diesen Fall berechnen zu k\"onnen, ist es offensichtlich, dass man eine geschlossene Kongruenzformel f\"ur elliptische Kurve \"uber einem endlichen K\"orper finden muss. Wir formulieren eine Frage von H. Esnault: Seien $p$ eine Primzahl, $\mathbb{F}_{p^n}$ eine endliche K\"orpererweiterung und $E/\mathbb{F}_p$ eine elliptische Kurve. Welche Elemente $a \in \mathbb{Z}/p^n$ in der Kongruenzformel
\[\#E(\mathbb{F}_{p^n}) \equiv a \mod p^n\]
f\"ur ein vorgegebenes $n$ aufgehen k\"onnen. Das Problem besteht darin, dass man i. Allg. nur schwer eine exakte geschlossene Formel f\"ur die Anzahl der rationalen Punkten von einer elliptischen Kurve bestimmen kann. Die Fragestellung von H. Esnault ist sehr eng verwandt mit einer offenen Frage von N. Koblitz (vgl. \cite[p. 164]{Kob88}), welche wiederum eine wichtige Bedeutung in den kryptographischen Anwendungen spielt, wenn wir $E$ und $p$ festhalten, dabei aber $n$ varieren lassen. Im Zusammenhang zur Frage von H. Esnault kann man nat\"urlich eine andere Vorgehensweise w\"ahlen. Man nimmt etwa die elliptischen Kurven von der Form $E_{A,B}: Y^2 = X^3 + AX + B$ \"uber $\mathbb{F}_p$ mit $\Delta(E_{A,B}) \neq 0$ und untersucht die Kongruenzformel
\[\#E_{A,B}(\mathbb{F}_{p^n}) \equiv a \mod p^n\]
f\"ur ein festes $n$, wobei $A$ und $B$ varieren. Es gilt
\[\#E_{A,B}(\mathbb{F}_{p^n}) = 1 + p^n + \sum_{x \in \mathbb{F}_{p^n}}\big ( \frac{x^3 + Ax + B}{p^n} \big ),\]
wobei 
\[\big ( \frac{}{p^n} \big ): \mathbb{F}_{p^n}^{\times} \rightarrow \mathbb{F}_{p^n}^{\times}/\mathbb{F}_{p^n}^{\times 2}\]
das Legendre-Symbol \"uber $\mathbb{F}_{p^n}$ ist. Da $n$ fest ist, kann man den folgenden einfachen Algorithmus anwenden, wobei wir die Effizienz au\ss er Acht lassen.
\begin{algorithm}
\caption{Zur Berechnung von $\#E_{A,B}(\mathbb{F}_{q}) \mod q, \quad q = p^n$}\label{Legendre}
\begin{algorithmic}[1]
\State $list \gets list()$
\BState \emph{loop}: \For{$A,B$ in $\mathbf{range}(\mathbb{F}_{p^n})}$
\State $\Delta(E_{A,B}) \gets 4A^3+27B^2$
\If{$\Delta(E_{A,B}) = 0$} 
\State $\mathbf{continue}$
\Else
\For{$x$ in $\mathbf{range}(\mathbb{F}_{p^n})$} 
\State $s \gets \big ( \frac{x^3 + Ax + B}{p^n} \big )$ \label{chi}
\State $list.\mathbf{append}(s)$ 
\EndFor
\EndIf
\EndFor
\Return $a \gets 1+\mathbf{sum}(list) \mod p^n$
\end{algorithmic}
\end{algorithm}
Es ist bekanntlich aus der elementaren Zahlentheorie, dass man im Schritt \ref{chi} von dem Algorithmus \ref{Legendre} effizient berechnen kann. Die Laufzeit von dieser Berechnung betr\"agt $\log^2(\max\{p^n,x^3+Ax+B\})$. Nat\"urlich kann man direkt den Schoof-Algorithmus, ein deterministischer polynomieller Algorithmus der Komplexit\"at $\mathcal{O}(q^8)$, sowie den SEA-Algorithmus der Komplexit\"at $\mathcal{O}(q^6)$, verwenden, aber wir werden uns nicht mehr weiter in dieses Thema vertiefen. An dieser Stelle m\"ochten wir darauf hinweisen, dass wir in dieser Arbeit nicht vorhaben, eine pr\"azise Antwort zu den oben genannten Problemen herauszufinden, sondern wir versuchen auf einer anderen Art und Weise der Fragestellung, gewisse asymptotische Aussagen zu treffen. Obwohl es in konkreten Beispielen, wie etwa wenn die Gleichung einer elliptischen Kurve vorgegeben ist, bekannte Algorithmen gibt, die uns die Anzahl der rationalen Punkte auf der Kurve liefern, behandeln wir die elliptischen Kurven ganz i. Allg. und zwar ganz unabh\"angig von deren Gleichungen. Nun wollen wir kurz beschreiben, wie die Arbeit eingeteilt wird. \"Uberall in dieser Arbeit verwenden wir bekannte Methoden aus der klassischen Funktionentheorie, analytischen Zahlentheorie und analytischen Kombinatorik. Im ersten Abschnitt berechnen wir einige asymptotische Ergebnisse \"uber endlichen K\"orpern, die sehr elementar sind. Erstaunlicherweise liefert die Methode mit erzeugenden Funktionen bzgl. der Cauchy'sche Integralformel aus der analytische Kombinatorik eine bessere Formulierung f\"ur eine obere Schranke von $\frac{\# E(\mathbb{F}_{p^n})}{\# E(\mathbb{F}_p)}$ als die Methode der Absch\"atzung von Exponentialsummen. Da wir uns nur f\"ur die Arithmetik der elliptischen Kurven interessieren, benutzen wir nur die klassische Verfassung von dem Primzahlsatz, die man wie in der aktuellen Forschung der modernen analytischen Zahlentheorie noch viel besser versch\"arfen kann. Au\ss erdem arbeiten wir \"uberall im lokalen Fall nur mit elliptischen Kurven, die \"uber einem Primk\"orper $\mathbb{F}_p$ definiert sind, um die Notationen zu vereinfachen, obwohl man ohne irgendwelche Ver\"anderung unsere Ergebnisse ganz allgemein auf einen endlichen K\"orper $\mathbb{F}_q$ mit $q = p^n$ \"ubertragen kann. Wir wollen dem Leser seine Aufmerksamkeit erwecken, dass einige elementare Ergebnisse sich in dieser Arbeit f\"ur supersingul\"are elliptische Kurven erheblich leichter als f\"ur gew\"ohnliche elliptischen Kurven zu zeigen lassen. Der bewegende Grund f\"ur uns, asymptotisches Verhalten von Punkten auf supersingul\"aren elliptischen Kurven bzgl. K\"orpererweiterungen zu untersuchen, ist die gro\ss e Rolle von supersingul\"aren elliptischen Kurven in der Post-Quanten-Kryptographie, wo man glaubt, dass das Problem mit der Berechnung einer Geheim-Isogenie von einer supersingul\"aren elliptischen Kurve auch f\"ur Quanten-Rechner hart sein sollte. Im zweiten Abschnitt erweitern wir die erzeugende Funktion der Folge $\{\#E(\mathbb{F}_{p^n})\}_{n \in \mathbb{N}}$ zu einer globalen Funktion, die die gemischten Eigenschaften von der Riemannschen Zetafunktion und der $L$-Funktion von einer elliptischen Kurve $E$ \"uber dem K\"orper der rationalen Zahlen besitzt. Wir erhalten hierbei eine asympotische Aussage f\"ur die globale erzeugende Funktion einer elliptischen Kurve \"uber $\mathbb{Q}$, die eine \"ahnliche Eigenschaft wie die von der $L$-Funktion $L(E;s)$ besitzt. Wir werden einige analytische Resultaten f\"ur die Hecke $L$-Funktionen der Modulformen vom Gewicht $2$ beweisen. Der Leser sollte sich vertraulich mit dem Grundbegriff von Modulformen und ganzen Funktionen wie etwa in \cite{FrBu06} machen. Unsere bevorstehende Arbeit f\"ur eine Anwendung in der Kryptographie und IT-Sicherheit basiert auf dieser rein mathematischen Arbeit mit der Motivation aus der Arbeit von M. Anshel und D. Goldfeld \cite{AG97}.  
\section{Einige elementare Ergebnisse}
\noindent Seien nun $E/\mathbb{F}_p$ eine elliptische Kurve \"uber einem Primk\"orper $\mathbb{F}_p$ und $\mathbb{F}_{p^n}/\mathbb{F}_p$ eine endliche K\"orpererweiterung. Das Theorem von A. Weil besagt, dass es eine algebraische Ganzzahl $\alpha \in \bar{\mathbb{Z}} \subset \mathbb{C}$ mit 
\[\#E(\mathbb{F}_p) = 1 + p - (\alpha + \bar{\alpha})\]
und
\[\#E(\mathbb{F}_{p^n}) = 1 + p^n - (\alpha^n + \bar{\alpha}^n)\]
gibt, wobei $\mathrm{Nm}(\alpha) = p$ und $\mathrm{Nm}: \bar{\mathbb{Z}} \rightarrow \mathbb{R}_{\geq 0}, z \mapsto z\bar{z}$ die Norm-Abbildung ist. Es gilt
\[\frac{\#E(\mathbb{F}_{p^n})}{\#E(\mathbb{F}_p)} = \mathrm{Nm}\big (\frac{\alpha^n-1}{\alpha-1}\big ) = \mathrm{Nm}\big (\sum_{j=0}^{n-1}\alpha^j \big ).\]
Betrachtet man $\alpha$ als eine komplexe Zahl, so l\"asst sich
\[\alpha = \sqrt{p}\cdot \exp(i\vartheta_E), \quad \vartheta_E \in [0,\pi]\]
schreiben. Seien $t = \alpha + \bar{\alpha}$ und $t_n = \alpha^n + \bar{\alpha}^n$. Da $\alpha^n + \bar{\alpha}^n = 2\sqrt{p^n}\cos(n\vartheta_E)$ und $\cos(n\vartheta_E) = T_n(\cos(\vartheta_E))$ das $n$-te Tchebycheff-Polynom der ersten Art ist, hat man eine rekursive Formel
\begin{equation}\label{e0}
t_{n+1} = t_n \cdot t - p\cdot t_{n-1}, \quad \forall n \geq 1, t_0 = 2, t_1 = t.
\end{equation}
Die dazugeh\"orige erzeugende Funktion lautet
\[A_E(z) = \frac{1}{1-t\cdot z + p\cdot z^2},\quad t_n = \frac{1}{2\pi i}\ointctrclockwise_{\mid z \mid = \rho} \frac{dz}{z^{n+1} \cdot (1-t \cdot z + p \cdot z^2)},\]
wobei es reicht, $\rho < \frac{1}{\sqrt{p}}$ anzufordern. Asymptotisch hat man aus der obigen rekursiven Formel
\[\#E(\mathbb{F}_{p^n}) \sim 1 + p^n - t^n, \quad n \to +\infty.\] 
Kongruenzformel modulo $p$ ist trivial. Es gilt
\[\#E(\mathbb{F}_{p^n}) \equiv 1-t^n \mod p, \quad \forall n \geq 1.\]
\begin{satz}\label{s1}
Seien $E/\mathbb{F}_p$ eine supersingul\"are elliptische Kurve mit $p \neq 2,3$ und $n \geq 1$ eine nat\"urliche Zahl. Dann hat man
\begin{equation}\label{e1}
\#E(\mathbb{F}_{p^n}) \equiv \begin{cases} 1 \mod p^n, \quad n \equiv 1 \mod 2 \\ 1 - 2p^k \mod p^n, \quad n = 2k,\, k \equiv 0 \mod 2 \\ 1 + 2p^k \mod p^n, \quad n = 2k, \, k \equiv 1 \mod 2 \end{cases}
\end{equation}
\end{satz}
\begin{proof}
Ist $E$ supersingul\"ar und $p \neq 2,3$, so ist $\#E(\mathbb{F}_p) = p + 1$, d.h $\cos(\vartheta_E) = 0$, also $\vartheta_E = \frac{\pi}{2}$. F\"ur eine endliche K\"orpererweiterung $\mathbb{F}_{p^n}/\mathbb{F}_p$ hat man 
\[\#E(\mathbb{F}_{p^n}) = 1 + p^n - 2\cdot \sqrt{p^n}\cos(n\vartheta_E).\]
Aus 
\begin{equation*}
\cos(n\frac{\pi}{2}) = \begin{cases} 0, \quad n \equiv 1 \mod 2 \\ 1, \quad n = 2k,\, k \equiv 0 \mod 2 \\ -1 \mod p^n, \quad n = 2k, \, k \equiv 1 \mod 2 \end{cases}
\end{equation*}
folgt die Behauptung.
\end{proof}
\noindent F\"ur gew\"ohnliche elliptische Kurven hat man
\begin{satz}\label{s2}
Sei $E/\mathbb{F}_p$ eine gew\"ohnliche elliptische Kurve. Es gilt
\[\#E(\mathbb{F}_{p^n}) \not \equiv 1 \mod p^n, \quad \forall n \geq 1.\]
\end{satz}
\begin{proof}
Angenommen, es g\"abe eine nat\"urliche Zahl $n \geq 1$, so dass 
\[\#E(\mathbb{F}_{p^n}) \equiv 1 \mod p^n.\]
Man hat 
\[\#E(\mathbb{F}_{p^n}) = 1 + p^n - \mathrm{Spur}(\mathbf{Frob}^n,\mathrm{End}(T_\ell(E))),\] 
wobei $T_{\ell}(E)$ der Tate-Modul von $E$ mit einer Primzahl $\ell \neq p$ und $\mathbf{Frob}$ der Frobenius-Endomorphismus von $E$ ist. Wir werden hier nicht \"uber den Tate-Modul sowie Frobenius-Endomorphismus diskutieren. Den interessierten Leser verweisen wir auf \cite{Sil09}. Nun muss 
\[p \mid \mathrm{Spur}(\mathbf{Frob}^n,\mathrm{End}(T_{\ell}(E))).\] 
Daher muss $E$ supersingul\"ar \"uber $\mathbb{F}_{p^n}$ sein. Dies f\"uhrt aber zum Widerspruch zur Tatsache, dass eine gew\"ohnliche elliptische Kurve bzgl. endlichen K\"orpererweiterungen erhaltend bleibt. Damit wurde die Behauptung bewiesen.
\end{proof}
\noindent Sei $\mathbb{P} \subset \mathbb{N}$ die Menge aller Primzahlen. Aus \ref{e1} vom Satz \ref{s1} oben sieht man, dass es f\"ur eine supersingul\"are elliptische Kurve $E/\mathbb{F}_p$ mit $p \neq 2,3$ gilt
\[\sum_{\mathclap{\substack{2 < \ell \leq x,\, \ell \in \mathbb{P}\\ \frac{\#E(\mathbb{F}_{p^\ell})}{\# E(\mathbb{F}_p)} \in \mathbb{P}}}}1 \quad = \quad \sum_{\mathclap{\substack{2 < \ell \leq x,\, \ell \in \mathbb{P} \\ \frac{p^{\ell}+1}{p+1} \in \mathbb{P}}}}1,\]
weil man in diesem Fall pr\"azis $\#E(\mathbb{F}_{p^{\ell}}) = 1 + p^{\ell}$ f\"ur eine ungerade Primzahl $\ell$ hat. Das in \cite[p. 164]{Kob88} diskutierte Problem ist in diesem Fall \"aquivalent zu dem Problem der Unendlichkeit der Wagstaff-Primzahlen der basis $p$. Wir werden uns in dieser Arbeit nicht mit diesem Problem auseinandersetzen. Zun\"achst ben\"otigen wir einen kleinen Hilfssatz, der in der analytischen Zahlentheorie oft zum Erfolg f\"uhrt.
\begin{hilfssatz}[Partielle Integration]\label{h1}
Seien $a_n \in \mathbb{C}$ und eine Folge $b_1 < b_2 < \cdots < b_n \to \infty$ reeller Zahlen mit dem Summator $A(x) = \sum_{b_n \leq x}a_n$. Sei $f: [b_1,x] \rightarrow \mathbb{C}$ eine stetige und st\"uckweise stetig differenzierbare Funktion. Dann hat man
\[\sum_{b_n \leq x}a_n\cdot f(b_n) = f(x)\cdot A(x) - \int_{b_1}^x A(t)\cdot f'(t)dt.\]
\end{hilfssatz}
\begin{proof}
Es gilt f\"ur Stieltjes-Integralen die partielle Integration
\[\sum_{b_n \leq x}a_nf(b_n) = \int_{b_1}^x f(t)dA(t) = f(x)\cdot A(x) - \int_{b_1}^xA(t)\cdot f'(t)dt. \]
Daraus folgt die Behauptung.
\end{proof} 
\noindent Mithilfe von Exponentialsumme k\"onnen wir absch\"atzen
\begin{satz}\label{s3}
Seien $E/\mathbb{F}_p$ eine elliptische Kurve und $\mathbb{F}_{p^n}/\mathbb{F}_p$ eine endliche K\"orpererweiterung. Sei $\vartheta_E \in [0,\pi]$ mit
\[1 + p - \#E(\mathbb{F}_p) = 2\sqrt{p}\cdot \cos(\vartheta_E).\]
F\"ur jedes $\varepsilon \in (0,1)$ gibt es eine Konstante $C(\varepsilon)$ unabh\"angig von $n$, so dass es gilt
\[\frac{\#E(\mathbb{F}_{p^n})}{\#E(\mathbb{F}_p)} \in [1,(p^{\frac{n-1}{2}}C(\varepsilon)((n-1)^{\varepsilon}+(n-2)^{\varepsilon})+2p^{\frac{n-1}{2}} - C(\varepsilon)(n-2)^{\varepsilon}-1)^2].\]
In Bezug auf $\vartheta_E$ hat man
\[\frac{\#E(\mathbb{F}_{p^n})}{\#E(\mathbb{F}_p)} \in [1,\big ( 1+\frac{2p^{\frac{n-1}{2}}}{|\sin(\vartheta_E/2)|}-\frac{1}{|\sin(\vartheta_E/2)|} \big )^2].\]
\end{satz}
\begin{proof}
F\"ur ein $u \in \mathbb{R}_{\geq 0}$ schreiben wir
\[S(u) = \sum_{j=0}^u \exp(2\pi i (j\cdot\frac{\vartheta_E}{2\pi})).\]
Nach dem Hilfssatz \ref{h1} hat man
\begin{equation*}
\begin{split}
\sum_{j=0}^{n-1}p^{j/2}\cdot \exp(2\pi i(j \cdot \frac{\vartheta_E}{2\pi})) = \int_{0}^{n-1}p^{u/2}dS(u) \\ = p^{\frac{n-1}{2}}\cdot \sum_{j=0}^{n-1}\exp(2\pi i(j\cdot \frac{\vartheta_E}{2\pi})) + \frac{1}{2}\log(p)\cdot \int_0^{n-1}S(u)p^{u/2}du.
\end{split}
\end{equation*}
Damit erhalten wir eine Ungleichung
\begin{equation*}
\mid \sum_{j=0}^{n-1}p^{j/2}\cdot \exp(2\pi i(j\cdot \frac{\vartheta_E}{2\pi})) \mid \leq p^{\frac{n-1}{2}}\cdot \mid S(n-1) \mid + \frac{1}{2}\log(p)\cdot \int_0^{n-1}|S(u)|p^{u/2}du.
\end{equation*}
Nach dem Weylschen Verfahren (vgl. \cite[\S I. 2]{Wal63}) gibt es zu jedem $\varepsilon \in (0,1)$ eine Konstante $C(\varepsilon)$, so dass es gilt
\begin{equation*}
\begin{split}
p^{\frac{n-1}{2}}\cdot \mid S(n-1) \mid + \frac{1}{2}\log(p)\cdot \int_0^{n-1}|S(u)|p^{u/2}du \leq p^{\frac{n-1}{2}}(C(\varepsilon)(n-1)^{\varepsilon}+1) + \\ \frac{1}{2}\log(p)\cdot\int_0^{n-1}p^{u/2}du + \frac{1}{2}\log(p)\cdot C(\varepsilon) \cdot \int_0^{n-1} [u]^{\varepsilon}p^{u/2}du
\end{split}
\end{equation*}
Die rechte Seite der Ungleichung ist gleich
\begin{equation*}
\begin{split}
p^{\frac{n-1}{2}}(C(\varepsilon)(n-1)^{\varepsilon}+1)+p^{\frac{n-1}{2}}-1 + C(\varepsilon)\cdot \sum_{k=1}^{n-1}\int_{k-1}^k [u]^{\varepsilon}p^{u/2}du = \\ p^{\frac{n-1}{2}}(C(\varepsilon)(n-1)^{\varepsilon}+1)+p^{\frac{n-1}{2}}-1 + \frac{1}{2}\log(p)\cdot C(\varepsilon)\cdot \sum_{k=1}^{n-1}(k-1)^{\varepsilon}\int_{k-1}^k p^{u/2}du = \\ p^{\frac{n-1}{2}}(C(\varepsilon)(n-1)^{\varepsilon}+1)+p^{\frac{n-1}{2}}-1 + C(\varepsilon)\cdot \sum_{k=1}^{n-1}(k-1)^{\varepsilon}(p^{k/2} - p^{\frac{k-1}{2}})
\end{split}
\end{equation*}
Wir haben offensichtlich
\[\sum_{k=1}^{n-1}(k-1)^{\varepsilon}(p^{k/2} - p^{\frac{k-1}{2}}) \leq (n-2)^{\varepsilon}\cdot \sum_{k=1}^{n-1}(p^{k/2}-p^{\frac{k-1}{2}}) = (n-2)^{\varepsilon}(p^{\frac{n-1}{2}}-1).\]
Am Ende des Tages bekommen wir schlie\ss lich die Absch\"atzung
\[\frac{\#E(\mathbb{F}_{p^n})}{\#E(\mathbb{F}_p)} \leq (p^{\frac{n-1}{2}}C(\varepsilon)((n-1)^{\varepsilon}+(n-2)^{\varepsilon})+2p^{\frac{n-1}{2}} - C(\varepsilon)(n-2)^{\varepsilon}-1)^2\]
Sei nun
\[S_1(u) = \sum_{j=1}^{u}\exp(2\pi i(j\cdot \frac{\vartheta_E}{2\pi})).\]
In Bezug auf $\vartheta_E$ kann man wegen dem Weylschen Verfahren nun absch\"atzen
\begin{equation*}
\begin{split}
\mid \sum_{j=0}^{n-1}\alpha^j \mid \leq 1 + \mid \int_1^{n-1}p^{u/2}dS_1(u) \mid \leq 1 + p^{u/2}\cdot|S_1(n-1)| + \frac{1}{2}\log(p)\cdot \int_{1}^{n-1}|S_1(u)|p^{u/2}du \\ \leq 1 + p^{(n-1)/2}\cdot \frac{1}{|\sin(\vartheta_E/2)|} + \frac{1}{2}\log(p)\cdot \frac{1}{|\sin(\vartheta_E/2)|}\cdot \int_1^{n-1}p^{u/2}du \\ = 1 + \frac{2p^{(n-1)/2}}{|\sin(\vartheta_E/2)|}-\frac{1}{|\sin(\vartheta_E/2)|}
\end{split}
\end{equation*}
und damit wurde der Satz bewiesen.
\end{proof}
\begin{hilfssatz}\label{h2}
Sei $E/\mathbb{F}_p$ eine elliptische Kurve. Seien $n_1, n_2 \in \mathbb{N}$ mit $n_1 \neq n_2$. Dann ist $\#E(\mathbb{F}_{p^{n_1}}) \neq \#E(\mathbb{F}_{p^{n_2}})$.
\end{hilfssatz}
\begin{proof}
Wir betrachten die Funktion $f: [2,\infty) \rightarrow \mathbb{R}_+$
\[f(x) = 1 + p^x - 2p^{x/2}\cdot \cos(x\vartheta_E), \quad \vartheta_E \in [0,\pi].\]
Ihre erste Ableitung ist
\begin{equation*}
\begin{split}
f^{\prime}(x) = p^x\cdot \log p + 2p^{x/2}(\vartheta_E\cdot \sin(x\vartheta_E) - \frac{1}{2} \cos(x\vartheta_E)\log p) > p^x\cdot \log p - 2p^{x/2}(\vartheta_E + \frac{1}{2} \log p) \\ > p^x \cdot \log p - 2p^{x/2}(\pi + 2 \log p) > p^{x} - p^{x/2} > 0.
\end{split}
\end{equation*}
Damit ist $f(x)$ streng monoton wachsend. Daraus folgt die Behauptung. 
\end{proof}
\noindent Sei im folgenden $\pi(x)$ die Primzahlfunktion
\[\pi(x) = \sum_{p \in \mathbb{P}, p \leq x}1 \]
\begin{bem}
{\rm
Der Satz \ref{s3} zusammen mit dem Hilfssatz \ref{h2} und dem Primzahlsatz mit Restglied (vgl. \cite[\S IV. 5]{Sch69}) liefern uns eine Absch\"atzung f\"ur eine elliptische Kurve $E/\mathbb{F}_p$
\begin{equation*}
\begin{split}
\sum_{\mathclap{\substack{\ell \leq x,\, \ell \in \mathbb{P} \\ \frac{\#E(\mathbb{F}_{p^{\ell}})}{\#E(\mathbb{F}_p)} \in \mathbb{P}}}}1 \leq \int_2^{(1+\frac{2p^{\frac{x-1}{2}}}{|\sin(\vartheta_E/2)|}-\frac{1}{|\sin(\vartheta_E/2)|})^2}\frac{du}{\log u} + \\ \underline{\mathcal{O}}(\big (1+\frac{2p^{\frac{x-1}{2}}}{|\sin(\vartheta_E/2)|}-\frac{1}{|\sin(\vartheta_E/2)|} \big )^2 \cdot \exp(-\frac{1}{200}\log^{1/2} \big (1+\frac{2p^{\frac{x-1}{2}}}{|\sin(\vartheta_E/2)|}-\frac{1}{|\sin(\vartheta_E/2)|}\big )^2)),
\end{split}
\end{equation*}
die viel schlechter als die triviale Absch\"atzung
\[\sum_{\mathclap{\substack{\ell \leq x,\, \ell \in \mathbb{P} \\ \frac{\#E(\mathbb{F}_{p^{\ell}})}{\#E(\mathbb{F}_p)} \in \mathbb{P}}}}1 \leq \int_2^{x}\frac{du}{\log u} + \underline{\mathcal{O}}_{x \to +\infty}(x\cdot \exp(-\frac{1}{200}\log^{1/2}x))\]
ist. Diese Beobachtung entspricht der offensichtlichen Tatsache, dass $\frac{\#E(\mathbb{F}_{p^n})}{\#E(\mathbb{F}_p)}$ sehr gross im Vergleich zu $n$ ist. 
}
\end{bem}
\noindent Der folgender Satz liefert eine etwas bessere Absch\"atzung als die Folgerung aus dem Satz \ref{s3}, aber selbstverst\"andlich immer noch schlechter als die triviale Absch\"atzung.
\begin{satz}\label{s4}
Sei $E/\mathbb{F}_p$ eine elliptische Kurve mit
\[\#E(\mathbb{F}_p) = 1 + p - t, \mid t \mid \leq 2\sqrt{p}, p \neq 2, 3\]
und $t < 0$. Dann gibt es ein $N_0 = N_0(\varepsilon)  \in \mathbb{N}$ zu jedem $\varepsilon \in (0,1)$, so dass es $\forall x \geq N_0$ gilt
\begin{equation*}
\begin{split}
\sum_{\mathclap{\substack{\ell \leq x,\, \ell \in \mathbb{P} \\ \frac{\#E(\mathbb{F}_{p^{\ell}})}{\#E(\mathbb{F}_p)} \in \mathbb{P}}}}1 \leq \pi((1+\varepsilon)(\frac{1}{p-2\sqrt{p}}+(p+2\sqrt{p})^{x-1}\cdot [x-1] + [x]((p+2\sqrt{p})^{[x]} - (5+2\sqrt{5})^x)) \\ - \pi((1-\varepsilon)\frac{126}{1+p})
\end{split}
\end{equation*}
\end{satz}
\begin{proof}
Sei $t \in \mathbb{Z}$ mit $\mid t \mid \leq 2\sqrt{p}$ und 
\[\#E(\mathbb{F}_p) = 1 + p - t.\]
Es gilt
\[\#E(\mathbb{F}_{p^n}) = 1 + p^n - \frac{1}{2\pi i}\ointctrclockwise_{\mid z \mid = \rho} \frac{dz}{z^{n+1}(1-tz + pz^2)}, \quad \forall n \geq 1, \rho < \frac{1}{\sqrt{p}}.\]
Asymptotisch hat man 
\[\#E(\mathbb{F}_{p^n}) \sim 1 + p^n - t^n, n \to +\infty.\]
$\forall \varepsilon > 0$, $\exists N_0 = N_0(\varepsilon) \in \mathbb{N}$, so dass es f\"ur $\forall n \geq N_0$ gelten 
\[\frac{\#E(\mathbb{F}_{p^n})}{\#E(\mathbb{F}_p)} \leq (1+\varepsilon)\frac{1+p^n-t^n}{1+p-t} \leq (1+\varepsilon)(\frac{1}{p+2\sqrt{p}} + \sum_{j=0}^{n-1}(p+2\sqrt{p})^j)\]
und f\"ur $n$ ungerade, $\varepsilon < 1$ und $t < 0$ 
\[\frac{\#E(\mathbb{F}_{p^n})}{\#E(\mathbb{F}_p)} \geq (1-\varepsilon)\frac{1+p^n-t^n}{1+p-t} \geq (1-\varepsilon)\frac{1+p^n}{1+p} \geq (1-\varepsilon)\frac{126}{1+p}.\]
Nach dem Hilfssatz \ref{h1} k\"onnen wir absch\"atzen
\begin{equation*}
\begin{split}
\sum_{j=0}^{x-1}(p+2\sqrt{p})^j = (p+2\sqrt{p})^{x-1}\cdot [x-1] - \log(p+2\sqrt{p})\cdot \int_0^x[u](p+2\sqrt{p})^udu \\ = (p+2\sqrt{p})^{x-1}\cdot [x-1] - \log(p+2\sqrt{p})\cdot \sum_{n=0}^{[x]-1}n \int_n^{n+1}(p+2\sqrt{p})^udu - [x]\int_{[x]}^x(p+2\sqrt{p})^udu \\ \leq (p+2\sqrt{p})^{x-1}\cdot [x-1] + [x]((p+2\sqrt{p})^{[x]}-(p+2\sqrt{p})^x) \\ \leq (p+2\sqrt{p})^{x-1}\cdot [x-1] + [x]((p+2\sqrt{p})^{[x]} - (5+2\sqrt{5})^x).
\end{split}
\end{equation*}
Daraus folgt die Behauptung.
\end{proof}
\noindent Zur\"uck zur Frage von H. Esnault wollen wir im folgenden absch\"atzen, welches $a \in [0,p^n-1]$, das in der Kongruenzformel
\[\#E(\mathbb{F}_{p^n}) \equiv a \mod p^n\]
f\"ur eine elliptische Kurve $E/\mathbb{F}_p$ aufgehen kann, wobei wir die Schwierigkeit der Frage vermindern, indem wir $n$ varieren. Dies bedeutet, dass wir die Folge $\{\#E(\mathbb{F}_{p^n})\}_{n \in \mathbb{N}}$ stets als eine arithmetische Funktion betrachten, die weder multiplikativ noch additiv ist. Dazu definieren wir eine Funktion f\"ur eine beliebige elliptische Kurve $E/\mathbb{F}_p$: 
\[[x]_{E,p,a} \quad \stackrel{defn}{=} \quad \sum_{\mathclap{\substack{2 \leq n \leq x \\ \#E(\mathbb{F}_{p^n}) \equiv a \mod p^n}}}1.\]
Offenbar $[x]_{E,p,a} \leq [x]-1$, wobei $[x]$ die Gau\ss sche Klammer ist. Nach der Behauptung \cite[Claim 3.1]{EsnXu09} hat man f\"ur jede elliptische Kurve $E/\mathbb{F}_p$
\begin{equation}\label{e2}
\lim_{x \to +\infty}\sum_{\mathclap{\substack{n \leq x \\ \#E(\mathbb{F}_{p^n}) \not\equiv 1 \mod p^n}}}1 \quad = +\infty
\end{equation}
Wegen dem Satz \ref{s1} hat man $\lim_{x \to +\infty}[x]_{E,p,1} = +\infty$, wenn $E/\mathbb{F}_p$ supersingul\"ar ist. Trotz der Gleichung \ref{e2} ist es nicht klar, ob es $\lim_{x \to +\infty}[x]_{E,p,a} \to +\infty$ f\"ur eine gew\"ohnliche elliptische Kurve $E/\mathbb{F}_p$ gilt und wenn, dann f\"ur welches $a$? Eine pr\"azise Antwort auf diese Frage k\"onnen wir nicht liefern, aber wir zeigen folgenden Satz, der eine asymptotische Aussage unter bestimmter Annahme trifft. 
\begin{satz}\label{s5}
Seien $E/\mathbb{F}_p$ eine elliptische Kurve mit
\[\#E(\mathbb{F}_p) = 1 + p - 2\cdot\sqrt{p}\cdot \cos(\vartheta_E), \quad \vartheta_E \in [0,\pi], p \neq 2,3\]
und $a \in \mathbb{N}$ eine nat\"urliche Zahl. Sei $x \in \mathbb{R}$ eine Zahl. Seien 
\[n_0 = \min_{n \in \mathbb{N} \cap [2,x]}\{n: \#E(\mathbb{F}_{p^n}) \equiv a \mod p^n\}, N = \max_{n \in \mathbb{N} \cap [2,x]}\{n: \#E(\mathbb{F}_{p^n}) \equiv a \mod p^n\}.\]
Dann ist entweder $a \leq 1 + 2 \cdot \sqrt{p^{n_0}}$ oder $p^N + 1 - 2\cdot \sqrt{p^N} \leq a \leq p^N$. Ferner, falls es gilt 
\[[x]_{E,p,a} \sim [x] - 1,\] 
dann gibt es zu jedem $\varepsilon \in (0,1)$ ein $\delta = \delta(\varepsilon) > 0$, so dass es f\"ur $x \to +\infty$ gelten
\[\frac{1}{2}(1-\varepsilon)p^{-\delta/2}-[x]+1-\vartheta_E(\varepsilon \cdot \frac{(\delta-1)^2}{2} + (1+\varepsilon)\cdot \frac{(x-1)^2}{2}) \leq \frac{a}{2}(1-\varepsilon)p^{-\delta/2}, \forall \varepsilon \in (0,1),\]
f\"ur $a < 1$ und
\begin{equation*}
\begin{split}
[x]-1 + \vartheta_E(\varepsilon \cdot \frac{(\delta-1)^2}{2} + (1+\varepsilon)\cdot \frac{(x-1)^2}{2})\\ \geq \frac{1-a}{2}(\delta(\frac{1}{p}-\frac{1}{p^{\delta/2}})+(1+\varepsilon)(\delta-1)p^{-\delta/2}+(1+\varepsilon)\frac{2}{\log p}p^{-\delta/2})
\end{split}
\end{equation*}
f\"ur $1 \leq a \leq 1+2\cdot \sqrt{p^{n_0}}$ und 
\[a = \underline{\mathcal{O}}_{x \to +\infty}(p^x), \quad a \geq p^N+1-2\cdot\sqrt{p^N}.\]
\end{satz}
\begin{proof}
Es gilt
\[\#E(\mathbb{F}_{p^n}) = 1 + p^n - 2\cdot \sqrt{p^n}\cdot \cos(n\vartheta_E).\]
Nach der Voraussetzung von dem Satz gibt es zu einer gewissen Anzahl von Zahlen $n \leq x$ ein $k_n$ derart, dass es gilt
\[k_n \cdot p^n + a = 1 - 2\cdot \sqrt{p^n}\cdot \cos(n\vartheta_E),\]
wenn $x$ hinreichend gross gew\"ahlt wird. Da $p \geq 5$ und $\mid \cos(n\vartheta_E) \mid \leq 1$, darf man daraus schlie\ss en, dass entweder $k_n = 0$ oder $k_n = -1$ ist, wobei es im ersten Fall $a \leq 1 + 2\cdot \sqrt{p^{n_0}}$ und im zweiten Fall $a \geq p^N + 1-2\cdot \sqrt{p^N}$ gilt. Im ersten Fall haben wir dann 
\begin{equation}\label{e3}
\sum_{n \leq x}{}^{\prime}\cos(n\vartheta_E) = (1-a)\cdot \sum_{n \leq x}{}^{\prime}\frac{1}{2\sqrt{p^n}},
\end{equation}
wobei es in der Summe $\sum{}^{\prime}$ \"uber allen $2 \leq n \leq x$ mit $\#E(\mathbb{F}_{p^n}) \equiv a \mod p^n$ aufsummiert wird. Nach dem Hilfssatz \ref{h1} hat man
\begin{equation*}
\begin{split}
\sum_{n \leq x}{}^{\prime}\cos(n\vartheta_E) = \cos(x\vartheta_E)\cdot [x]_{E,p,a} + \vartheta_E \cdot \int_2^x [u]_{E,p,a}\cdot \sin(u\vartheta_E)du
\end{split}
\end{equation*}
$\forall \varepsilon > 0$, $\exists \delta = \delta(\varepsilon)$, so dass es $\forall u > \delta$ gilt $[u]_{E,p,a} \leq (1+\varepsilon)\cdot ([u]-1)$. Wir sch\"atzen ab
\[\vert \int_2^{\delta} [u]_{E,p,a} \cdot \sin(u\vartheta_E)du \vert \leq \int_{2}^{\delta}[u]_{E,p,a}du \leq \int_2^{\delta}([u]-1)du \leq \frac{(\delta-1)^2}{2}.\]
Damit haben wir
\begin{equation*}
\begin{split}
\sum_{n \leq x}{}^{\prime}\cos(n\vartheta_E) \leq \cos(x\vartheta_E) [x]_{E,p,a} + \vartheta_E \frac{(\delta-1)^2}{2} + \vartheta_E  \int_{\delta}^x[u]_{E,p,a} \sin(u\vartheta_E)du \\ \leq \cos(x\vartheta_E) [x]_{E,p,a} + \vartheta_E \frac{(\delta-1)^2}{2}+(1+\varepsilon) \vartheta_E  \int_{\delta}^x ([u]-1)du \\ \leq \cos(x\vartheta_E) [x]_{E,p,a} + \vartheta_E  \frac{(\delta-1)^2}{2} + (1+\varepsilon) \vartheta_E  \int_{\delta}^x (u-1)du \\ = \cos(x\vartheta_E) [x]_{E,p,a} + \vartheta_E \frac{(\delta-1)^2}{2} + (1+\varepsilon)\cdot \vartheta_E\cdot \frac{(x-1)^2}{2}-(1+\varepsilon)\cdot \vartheta_E \cdot \frac{(\delta-1)^2}{2} \\ \leq [x]-1 + \vartheta_E(\varepsilon \cdot \frac{(\delta-1)^2}{2} + (1+\varepsilon)\cdot \frac{(x-1)^2}{2}).
\end{split}
\end{equation*}
Nun sch\"atzen wir die rechte Seite der Gleichung \ref{e3} ab. Nach dem Hilfssatz \ref{h1} haben wir
\begin{equation*}
\begin{split}
\frac{1-a}{2}\sum_{n \leq x}{}^{\prime}\frac{1}{\sqrt{p^n}} = \frac{1-a}{2}\cdot (\frac{1}{p^{x/2}}\cdot [x]_{E,p,a} + \frac{\log p}{2}\cdot \int_2^x[u]_{E,p,a}\cdot p^{-u/2}du)
\end{split}
\end{equation*}
$\forall \varepsilon > 0$, $\exists \delta = \delta(\varepsilon) > 0$, so dass es $\forall u > \delta$ gilt $[u]_{E,p,a} \geq (1-\varepsilon)\cdot([u]-1)$. Damit erhalten wir f\"ur alle $\varepsilon$ mit $0 < \varepsilon < 1$ und $a \leq 1$
\begin{equation*}
\begin{split}
\frac{1-a}{2}\sum_{n \leq x}{}^{\prime}\frac{1}{\sqrt{p^n}} \geq \\ \frac{1-a}{2}(\frac{1}{p^{x/2}}(1-\varepsilon)([x]-1)+\frac{\log p}{2}\int_2^{\delta}[u]_{E,p,a}p^{-u/2}du + (1-\varepsilon)\frac{\log p}{2}\int_{\delta}^x([u]-1)p^{-u/2}du) \\ \geq \frac{1-a}{2}(\frac{1}{p^{x/2}}(1-\varepsilon)([x]-1) - (1-\varepsilon)\frac{\log p}{2}\cdot \int_{\delta}^x p^{-u/2}du) = \\ \frac{1-a}{2}(\frac{1}{p^{x/2}}(1-\varepsilon)([x]-2)+(1-\varepsilon)p^{-\delta/2}).
\end{split}
\end{equation*}
Daraus folgt f\"ur $x \to +\infty$, $a \leq 1$
\[\frac{1}{2}(1-\varepsilon)p^{-\delta/2}-[x]+1-\vartheta_E(\varepsilon \cdot \frac{(\delta-1)^2}{2} + (1+\varepsilon)\cdot \frac{(x-1)^2}{2}) \leq \frac{a}{2}(1-\varepsilon)p^{-\delta/2}, \forall \varepsilon \in (0,1).\]
Nun 
\begin{equation*}
\begin{split}
\sum_{n \leq x}{}^{\prime}\frac{1}{\sqrt{p^n}} = \frac{1}{p^{x/2}}[x]_{E,p,a}+\frac{\log p}{2}\int_2^x[u]_{E,p,a} \cdot p^{-u/2}du \\ \leq \frac{1}{p^{x/2}}(1+\varepsilon)([x]-1) + \frac{\log p}{2}\int_{2}^{\delta}[u]_{E,p,a}\cdot p^{-u/2}du + \frac{\log p}{2}(1+\varepsilon)\int_{\delta}^x ([u]-1)p^{-u/2}du \\ \leq \frac{1}{p^{x/2}}(1+\varepsilon)([x]-1) + \delta \cdot (\frac{1}{p}-p^{-\delta/2}) + \frac{\log p}{2}(1+\varepsilon)\int_{\delta}^x(u-1)p^{-u/2}du \\ = \frac{1}{p^{x/2}}(1+\varepsilon)([x]-1) + \delta \cdot (\frac{1}{p}-p^{-\delta/2}) + (1+\varepsilon)(\delta-1)p^{-\delta/2} - (1+\varepsilon)(\delta-1)(x-1)p^{-x/2} \\ - (1+\varepsilon)\frac{2}{\log p}p^{-x/2} + (1+\varepsilon)\frac{2}{\log p}p^{-\delta/2}
\end{split}
\end{equation*}
Daraus folgt f\"ur $x \to +\infty$, $a \geq 1$
\begin{equation*}
\begin{split}
[x]-1 + \vartheta_E(\varepsilon \cdot \frac{(\delta-1)^2}{2} + (1+\varepsilon)\cdot \frac{(x-1)^2}{2})\\ \geq \frac{1-a}{2}(\delta(\frac{1}{p}-\frac{1}{p^{\delta/2}})+(1+\varepsilon)(\delta-1)p^{-\delta/2}+(1+\varepsilon)\frac{2}{\log p}p^{-\delta/2})
\end{split}
\end{equation*}
Im zweiten Fall haben wir dann
\begin{equation}\label{e4}
\sum_{n \leq x}{}^{\prime}\cos(n\vartheta_E) = \sum_{n \leq x}{}^{\prime}\frac{p^n+1-a}{2\sqrt{p^n}} = \frac{1-a}{2}\sum_{n \leq x}{}^{\prime}\frac{1}{\sqrt{p^n}} + \frac{1}{2}\sum_{n \leq x}{}^{\prime}\sqrt{p^n}
\end{equation}
mit $a \geq p^n+1-2\sqrt{p^n}$. Nach der Voraussetzung von dem Satz gibt es zu jedem $\varepsilon > 0$ ein $\delta = \delta(\varepsilon) > 0$, so dass es f\"ur alle $u > \delta$ gilt
\begin{equation*}
(1-\varepsilon)\cdot([u]-1) \leq [u]_{E,p,a} \leq (1+\varepsilon)\cdot([u]-1).
\end{equation*}
Es gelten nach dem Hilfssatz \ref{h1}
\begin{equation*}
\begin{split}
\sum_{n\leq x}{}^{\prime}\cos(n\vartheta_E) = \cos(x\vartheta_E)[x]_{E,p,a}+\vartheta_E\int_2^x[u]_{E,p,a}\sin(u\vartheta_E)du \\ \geq -(1+\varepsilon)([x]-1) - \vartheta_E\frac{(\delta-1)^2}{2} + \vartheta_E(1-\varepsilon)\int_{\delta}^x([u]-1)\sin(u\vartheta_E)du \\ \geq -(1+\varepsilon)([x]-1) - \vartheta_E\frac{(\delta-1)^2}{2} - \vartheta_E(1-\varepsilon)\int_{\delta}^x([u]-1)du \\ \geq -(1+\varepsilon)([x]-1) - \vartheta_E\frac{(\delta-1)^2}{2} - \vartheta_E(1-\varepsilon)\int_{\delta}^x (u-1)du \\ \geq -(1+\varepsilon)([x]-1) - \vartheta_E\frac{(\delta-1)^2}{2} - \vartheta_E(1-\varepsilon)\frac{(x-1)^2}{2} + \vartheta_E(1-\varepsilon)\frac{(\delta-1)^2}{2} \\ = -(1+\varepsilon)([x]-1) - \vartheta_E\varepsilon\frac{(\delta-1)^2}{2}-\vartheta_E\frac{(x-1)^2}{2}, \forall \varepsilon \in (0,1)
\end{split}
\end{equation*}
Wir sch\"atzen die rechte Seite der Gleichung \ref{e4} f\"ur $0 < \varepsilon < 1$ ab
\begin{equation*}
\begin{split}
\sum_{n \leq x}{}^{\prime} \frac{p^n+1-a}{2\sqrt{p^n}} \leq \frac{1-a}{2}\sum_{n\leq x}{}^{\prime}\frac{1}{\sqrt{p^n}} + p^{x/2}[x]_{E,p,a} \\ \leq p^{x/2}(1+\varepsilon)([x]-1) + \frac{1-a}{2}(\frac{1}{p^{x/2}}(1-\varepsilon)([x]-2)+(1-\varepsilon)p^{-\delta/2})
\end{split}
\end{equation*}
In diesem Fall liefert unsere Methode die trivial asymptotische Absch\"atzung $a = \underline{\mathcal{O}}_{x\to +\infty}(p^x)$. Damit wurde der Satz bewiesen.
\end{proof}
\noindent Eine Funktion $L: [a,+\infty] \rightarrow \mathbb{R}$ mit $a \geq 0$ hei\ss t von langsamem Wachstum, wenn $L \geq 0$ und stetig ist, und wenn es f\"ur alle $c \in (0,+\infty)$ gilt
\[\lim_{x \to +\infty}\frac{L(cx)}{L(x)} = 1.\] 
Sei im folgenden $\Gamma(s)$ die Gammafunktion
\[\Gamma(s) = \int_0^{+\infty}t^{s-1}e^{-t}dt, \quad \mathfrak{R}(s) > 0.\]
Es gilt nach partieller Integration $\Gamma(n+1) = n!, \forall n \in \mathbb{N}_+$. 
\begin{hilfssatz}[Hardy-Littlewood-Karamata]\label{h3}
Sei $A(\cdot)$ eine $\mathbb{R}$-, monoton nicht-fallende Funktion auf $[0,+\infty)$ mit $A(0) = 0$. Angenommen, das Laplace-Integral
\[f(\sigma) = \sigma \int_0^{+\infty}A(t)\cdot \exp(-\sigma t)dt < \infty, \quad \forall \sigma > 0.\]
Ferner, seien $L$ eine Funktion von langsamen Wachstum und $\alpha$ eine nicht-negative Konstante mit
\[\lim_{\sigma \to 0+}f(\sigma)\cdot \sigma^{-\alpha}\cdot \frac{1}{L(1/\sigma)} = A.\]
Dann hat man die asympotische Formel
\[A(x) \sim \frac{A}{\Gamma(\alpha+1)}\cdot x^{\alpha} \cdot L(x), \quad x \to +\infty.\]
\end{hilfssatz}
\noindent Der Hilfssatz \ref{h3} ist bekannt unter dem Name Taubersatz von Hardy-Littlewood-Karamata (vgl. \cite[\S V.4]{Sch69}). F\"ur eine in $\mid \sigma \mid < 1$ konvergierte Potenzreihe
\[f(\sigma) = \sum_{n \geq 0}a_n\sigma^n, \quad a_n \geq 0\]
kann man umformen
\[f(\sigma) = \log \frac{1}{\sigma} \cdot \int_0^{+\infty} \big (\sum_{n < u}a_n \big )\exp(-u\log \frac{1}{\sigma})du.\]
Der Hilfssatz \ref{h3} liefert dann eine asympotisch Formel
\[\sum_{n < x}a_n \sim \frac{A}{\Gamma(\alpha+1)}\cdot x^{\alpha}\cdot L(x), \quad x \to +\infty,\]
falls
\[f(\sigma) \sim A(1-\sigma)^{-\alpha}\cdot L((1-\sigma)^{-1}), \quad \sigma \to 1^{-}.\]
Mit dem Hilfssatz \ref{h3} und mithilfe von erzeugenden Funktionen aus der analytischen Kombinatorik kann man folgenden Satz rein analytisch zeigen, ohne den Satz von Hasse f\"ur elliptische Kurven anzuwenden. Vorausgesetzt, man kennt die rekursive Formel zuvor. 
\begin{satz}\label{s6}
Sei $E/\mathbb{F}_p$ eine elliptische Kurve. Dann hat man asymptotisch
\[\sum_{n \leq x}\frac{\#E(\mathbb{F}_{p^n})}{p^n} \sim x, \quad x \to +\infty\]
und
\[\frac{\sum_{n \leq x}\#E(\mathbb{F}_{p^n})\cdot p^{-n}}{\pi(x)} \sim \log x, \quad x \to +\infty.\]
\end{satz}
\begin{proof}
Die zweite asymptotische Formel folgt aus der ersten Formel und dem Primzahlsatz
\[\pi(x) = \int_2^x \frac{du}{\log u} + \underline{\mathcal{O}}_{x \to +\infty}(x\cdot \exp(-\frac{1}{200}\log^{1/2} x)) \sim \frac{x}{\log x}, \quad x \to +\infty.\]
Wir zeigen nun die erste Aussage. Sei $\{e_n\}_{n \in \mathbb{N}}$ die Folge $\{\#E(\mathbb{F}_{p^n})\}_{n \in \mathbb{N}}$ mit $e_n = 0$ f\"ur $n \leq 0$. Es gilt
\[e_n = 1 + p^n - t_n, \quad \forall n \geq 1\]
mit
\[t_{n+1} = t_n \cdot t - p\cdot t_{n-1}, \quad \forall n \geq 1, t_0 = 2, t_1 = t.\]
Sei $B_E(\sigma) = \sum_{n=0}^{+\infty}e_n\sigma^n$ die erzeugende Funktion von der Folge $\{e_n\}_{n\in \mathbb{N}}$, d.h
\[e_n = \frac{1}{2 \pi i} \ointctrclockwise_{\mid z \mid = \rho} \frac{B_E(z)}{z^{n+1}} dz\]
f\"ur ein geeignetes $\rho$. Die Potenzreihe $B_E(\sigma)$ konvergiert in dem Bereich $\mid \sigma \mid < \frac{1}{p}$. Es gilt
\[e_n = t\cdot e_{n-1} - p\cdot e_{n-2} + (1+p^{n-1})\cdot e_1, \quad \forall n \geq 2.\]
Aus der rekursiven Formel leitet man die dazugeh\"orige erzeugende Funktion ab
\[B_E(\sigma) = \frac{e_1 \sigma}{1-t\sigma+p\sigma^2}\big ( \frac{1}{1-\sigma} + \frac{p\sigma}{1-p\sigma}\big ). \]
Sei $f(v) = \sum_{n=0}^{+\infty}e_n{}^{\prime}\cdot v^n$ mit $e_n{}^{\prime} = \frac{e_n}{p^n}$ und $v = p\cdot \sigma$. Es ist offenbar, dass $f(v)$ dann in dem Bereich $\mid v \mid < 1$ konvergiert.
Es gilt
\[f(v) \sim \frac{1}{1-v}, \quad v \to 1^{-}.\]
Der Taubersatz von Hardy-Littlewood-Karamata \ref{h3} liefert uns dann die gew\"unschte Asymptotik
\[\sum_{n \leq x}\#E(\mathbb{F}_{p^n})\cdot p^{-n} \sim x, \quad x \to +\infty.\]
Damit wurde der Satz bewiesen.
\end{proof}
\noindent Wir wollen im folgenden eine bessere Formulierung f\"ur die Absch\"atzung in dem Satz \ref{s3} zeigen. Es ist nicht nur, dass die Absch\"atzung in diesem Satz noch sch\"oner formuliert wird, sondern l\"asst sich die positive Konstante auch dabei sehr einfach bestimmen.
\begin{satz}\label{s7}
Sei $E/\mathbb{F}_p$ eine elliptische Kurve mit 
\[\#E(\mathbb{F}_p) = 1 + p - t, \quad \mid t \mid \leq 2\sqrt{p}.\]
$\forall \varepsilon > 0$, $\exists \, C = C(p,\varepsilon,t)$ eine positive Konstante, so dass
\[\frac{\# E(\mathbb{F}_{p^n})}{\# E(\mathbb{F}_p)} \leq C\cdot (p^{1/2+\varepsilon})^{n-1}.\]
\end{satz}
\begin{proof}
Wie in dem Beweis des Satzes \ref{s6} hat man
\[\frac{\#E(\mathbb{F}_{p^n})}{\#E(\mathbb{F}_p)} = \frac{1}{2 \pi i} \ointctrclockwise_{\mid z \mid = \rho} \frac{1}{z^n(1-tz+pz^2)}\big ( \frac{1}{1-z} + \frac{pz}{1-pz} \big )dz, \quad \rho < \frac{1}{\sqrt{p}}.\]
Wir haben
\[\mid 1-tz+pz^2\mid \geq \mid 1-\mid t \mid\rho - p\cdot \rho^2 \mid, \mid 1-z \mid \geq 1 - \rho, \mid 1 - pz \mid \geq \mid 1 - p\cdot \rho \mid.\]
Es gilt
\begin{equation*}
\begin{split}
\big \vert \frac{1}{2 \pi i} \ointctrclockwise_{\mid z \mid = \rho} \frac{1}{z^n(1-tz+pz^2)}\big ( \frac{1}{1-z} + \frac{pz}{1-pz} \big )dz \big \vert \leq \\ \frac{1}{2\pi} \int_0^{2\pi}\big\vert  \frac{1}{z^n(1-tz+pz^2)}\big ( \frac{1}{1-z} + \frac{pz}{1-pz} \big ) \big\vert \rho d\theta \leq \frac{1}{\rho^{n-1}}\cdot \frac{1}{\mid 1-\mid t\mid \cdot \rho - p\rho^2 \mid}(\frac{1}{1-\rho}+\frac{\sqrt{p}}{\mid 1-p\rho \mid}).
\end{split}
\end{equation*}
Da $\rho < \frac{1}{\sqrt{p}}$, so ist $\rho \leq \frac{1}{p^{1/2+\varepsilon}}, \forall \varepsilon > 0$. W\"ahle $\rho = \frac{1}{p^{1/2+\varepsilon}}$. Dann hat man
\[\frac{1}{\rho^{n-1}}\cdot \frac{1}{\mid 1-\mid t\mid \cdot \rho - p\rho^2 \mid}(\frac{1}{1-\rho}+\frac{\sqrt{p}}{\mid 1-p\rho \mid}) = C(p,\varepsilon,t)\cdot (p^{1/2+\varepsilon})^{n-1},\]
wobei
\[C(p,\varepsilon,t)=\frac{\mid 1-\frac{1}{p^{-1/2+\varepsilon}}\mid + \sqrt{p}(1-\frac{1}{p^{1/2+\varepsilon}})}{\mid 1-\mid t\mid \cdot \frac{1}{p^{1/2+\varepsilon}}-\frac{1}{p^{2\varepsilon}}\mid \cdot (1-\frac{1}{p^{1/2+\varepsilon}}) \cdot \big \vert 1-\frac{1}{p^{-1/2+\varepsilon}} \big \vert}.\]
Daraus folgt die Behauptung.
\end{proof}
\begin{folgerung}\label{f1}
Sei $E/\mathbb{F}_p$ eine elliptische Kurve mit $\#E(\mathbb{F}_p) = 1 + p - t$. Zu jedem $\varepsilon > 0$ gibt es eine positive Konstante $C = C(p,\varepsilon,t)$, so dass
\begin{equation*}
\begin{split}
\sum_{\mathclap{\substack{\ell \leq x, \, \ell \in \mathbb{P} \\ \frac{\# E(\mathbb{F}_{p^{\ell}})}{\#E(\mathbb{F}_p)}\in \mathbb{P}}}}1 \leq \int_2^{C\cdot (p^{1/2+\varepsilon})^{x-1}}\frac{du}{\log u} + \\ \underline{\mathcal{O}}_{x \to +\infty}(p^{(1/2+\varepsilon)(x-1)}\exp(-\frac{1}{200}\sqrt{(\frac{1}{2}+\varepsilon)(x-1)\log p + \log C}))
\end{split}
\end{equation*}
\end{folgerung}
\begin{proof}
Die Behauptung folgt aus dem Satz \ref{s7} und dem Primzahlsatz mit Restglied. 
\end{proof}
\begin{satz}\label{s8}
Seien $E/\mathbb{F}_p$ eine elliptische Kurve und $x \in \mathbb{R}_+$ eine Zahl. Es gibt dann ein $\theta = \theta(x) \in [-\infty,1]$, so dass 
\[\sum_{\mathclap{\substack{\ell \leq x,\, \ell \in \mathbb{P} \\ \frac{\# E(\mathbb{F}_{p^{\ell}})}{\#E(\mathbb{F}_p)}\in \mathbb{P}}}}1 = \pi(x^{\theta}).\]
Insbesondere, 
\[\sum_{\mathclap{\substack{\ell \leq x,\, \ell \in \mathbb{P} \\ \frac{\# E(\mathbb{F}_{p^{\ell}})}{\#E(\mathbb{F}_p)}\in \mathbb{P}}}}1 = \int_2^{x^{\theta}}\frac{du}{\log u} + \underline{\mathcal{O}}_{x \to +\infty}(x^{\theta}\exp(-\frac{1}{200}\sqrt{\log x^{\theta}}))\]
und
\[\lim_{x \to +\infty} \quad \sum_{\mathclap{\substack{\ell \leq x,\, \ell \in \mathbb{P} \\ \frac{\# E(\mathbb{F}_{p^{\ell}})}{\#E(\mathbb{F}_p)}\in \mathbb{P}}}}1 = +\infty,\] 
falls $\theta > 0$.
\end{satz}
\begin{proof}
Die Existenz von $\theta$ ist klar, n\"amlich
\[\theta = \log_x(\max_{\ell \leq x, \, \ell \in \mathbb{P}}\{\ell: \frac{\#E(\mathbb{F}_{p^{\ell}})}{\#E(\mathbb{F}_p)} \in \mathbb{P}\}).\]
Es gilt daher $\max_{\ell \leq x, \, \ell \in \mathbb{P}}\{\ell: \frac{\#E(\mathbb{F}_{p^{\ell}})}{\#E(\mathbb{F}_p)} \in \mathbb{P}\} = [x^{\theta}]$. So muss 
\[\sum_{\mathclap{\substack{\ell \leq x,\, \ell \in \mathbb{P} \\ \frac{\# E(\mathbb{F}_{p^{\ell}})}{\#E(\mathbb{F}_p)}\in \mathbb{P}}}}1 = \pi(x^{\theta}).\]
Falls $\theta > 0$, hat man die Asymptotik aus dem Primzahlsatz
\[\pi(x^{\theta}) \sim \frac{x^{\theta}}{\theta\log x}, \quad x \to +\infty.\]
F\"ur jedes $\varepsilon>0$ hat man $\log(x) = \mathcal{O}_{\varepsilon}(x^{\varepsilon}), x \to +\infty$. Damit wurde der Satz bewiesen.
\end{proof}
\begin{satz}\label{s9}
Seien $E/\mathbb{F}_p$ eine elliptische Kurve \"uber einem Primk\"orper mit $p \neq 2,3$ und $B>0$ eine Zahl vorgegeben. Dann hat man
\begin{enumerate}
\item \label{(1)} $\forall x \geq e^{\sqrt[B]{2}}$, $\exists \, C=C(B)$, so dass 
\[\int_2^{x}\frac{du}{\log u} - C\cdot x \cdot e^{-\frac{1}{200}\log^{1/2}x} \leq \quad \sum_{\mathclap{\substack{\ell \in \mathbb{P}, \ell \leq x \\ \#E(\mathbb{F}_{p^{\ell}}) \equiv 1 \mod p^\ell}}}1 \quad \leq \int_2^x\frac{du}{\log u} + C\cdot x \cdot e^{-\frac{1}{200}\log^{1/2}x},\]
falls $E$ supersingul\"ar ist. 
\item \label{(2)} Falls $E$ gew\"ohnlich ist, dann
\[\sum_{\mathclap{\substack{\ell \in \mathbb{P}, \ell \leq x \\ \#E(\mathbb{F}_{p^{\ell}}) \not\equiv 1 \mod p^\ell}}}1 \quad = \quad \int_2^x \frac{du}{\log u} + \underline{\mathcal{O}}_{x \to +\infty}(x\cdot \exp(-\frac{1}{200}\log^{1/2}x)), \quad x \to +\infty \]
\end{enumerate} 
\end{satz}
\begin{proof}
Der Satz folgt automatisch aus den S\"atzen \ref{s1}, \ref{s2} und dem Primzahlsatz von Siegel-Walfisz (vgl. \cite[\S VI]{Sch69})
\end{proof}
\section{Globale erzeugende Funktion}
\noindent Das Ziel in diesem Abschnitt ist, die vorherige Diskussion \"uber die erzeugende Funktion von einer elliptischen Kurve \"uber einem endlichen K\"orper auf globalen Zustand zu erweitern. Wir beschr\"anken uns nur auf dem K\"orper der rationalen Zahlen $\mathbb{Q}$. F\"ur eine elliptische Kurve $E/\mathbb{F}_p$ mit $\#E(\mathbb{F}_p) = 1+p-t, \mid t \mid \leq 2\sqrt{p}$ ist es bekanntlich (vgl. \cite{Sil09}), dass ihre Zetafunktion
\[\zeta_{E/\mathbb{F}_p}(s) \stackrel{defn}{=} \frac{1-t\cdot p^{-s} + p^{1-2s}}{(1-p^{-s})(1-p^{1-s})}\]
die Funktionalgleichung
\[\zeta_{E/\mathbb{F}_p}(1-s) = \zeta_{E/\mathbb{F}_p}(s).\]
Die Riemannsche Vermutung f\"ur $E$ besagt, dass die Nullstellen $s$ von $\zeta_{E/\mathbb{F}_p}(s)$ den reellen Teil $\mathfrak{R}(s) = \frac{1}{2}$ haben. Die erzeugende Funktion
\[B_E(z) = \frac{\#E(\mathbb{F}_p)\cdot z}{(1-t\cdot z + p\cdot z^2)}\cdot \frac{1- p\cdot z^2}{(1-z)\cdot (1-p\cdot z)}\]
von der Folge $\{\#E(\mathbb{F}_{p^n})\}_{n \in \mathbb{N}}$ erf\"ullt die Funktionalgleichung
\[B_E(z) = -B_E \big (\frac{1}{p\cdot z} \big ),\]
wie man leicht nachpr\"ufen kann. Setzt man $z = p^{-s}$ mit $s = \sigma + i\cdot t \in \mathbb{C}$, so erh\"alt man die Funktionalgleichung
\begin{equation}\label{e5}
B_{E,p}(s) = -B_{E,p}(1-s)
\end{equation}
mit
\[B_{E,p}(s) = \frac{\# E(\mathbb{F}_p)p^{-s}}{1-tp^{-s}+p^{1-2s}}\cdot \frac{1-p^{1-2s}}{(1-p^{-s})(1-p^{1-s})}.\]
Die Menge der Nullstellen von $B_{E,p}(s)$ ist
\[\{ s \in \mathbb{C}: \mathfrak{R}(s) = \frac{1}{2},\, \mathfrak{Im}(s) \in \frac{\pi}{\log p}\cdot \mathbb{Z}\}\setminus \{s \in \mathbb{C}: \zeta_{E/\mathbb{F}_p}(s) = 0\}.\]
Im folgenden wollen wir die globale Situation betrachten, wo man $E/\mathbb{Q}$ eine \"uber $\mathbb{Q}$ definierte elliptische Kurve festh\"alt und dabei $p$ varieren l\"asst. Sei
\[\zeta(s) \stackrel{defn}{=} \sum_{n=1}^{+\infty}\frac{1}{n^s} = \prod_{p \in \mathbb{P}}\frac{1}{1-p^{-s}}, \quad \sigma > 1\]
die Riemannsche Zetafunktion, wobei die zweite Gleichheit die Eulersche Produktformel ist. Es ist bekanntlich, dass $\zeta(s)$ sich auf $\mathbb{C}-\{1\}$ analytisch fortsetzen l\"asst und einen einfachen Pol an der Stelle $s=1$ mit dem Residuum $\mathrm{Res}_{s=1}(\zeta(s)) = 1$ besitzt. Dies kann man leicht mit der Technik der partiellen Integration in dem Hilfssatz \ref{h1} zusammen mit dem Lebesgues Satz von der dominierten Konvergenz zeigen.  Ferner, $\zeta(s)$ erf\"ullt die Funktionalgleichung
\[\zeta(1-s) = 2(2\pi)^{-s}\Gamma(s)\cos(\frac{\pi s}{2})\zeta(s).\]
Die ber\"uhmte Riemannsche Vermutung behauptet, dass die nicht-trivialen Nullstellen von der Riemannschen Zetafunktion sich auf der kritischen Gerade $\sigma = \frac{1}{2}$ konzentrieren. Die M\"obiussche Funktion $\mu: \mathbb{N} \rightarrow \{0,\pm 1\}$ wird durch
\begin{equation*}
\mu(n) = \begin{cases}1, \quad n = 1 \\ 0, \quad \exists \, p \in \mathbb{P}, p^2 \, \mid \, n \\ (-1)^k, \quad n = p_1 \cdots p_k, \, \forall i,j \in \{1,\cdots,k\} \, p_i \neq p_j \end{cases}
\end{equation*}
definiert. F\"ur $\sigma > 1$ hat man
\[\frac{1}{\zeta(s)} = \sum_{n=1}^{+\infty}\frac{\mu(n)}{n^s}.\]
F\"ur eine elliptische Kurve $E/\mathbb{Q}$ schreiben wir $E_p$ f\"ur die Reduktion modulo $p$ von $E$. Es gilt f\"ur die guten Reduktionen
\[\#E_p(\mathbb{F}_p) = 1 + p - t_p.\]
Die $L$-Funktion auf $E$ ist per Definition durch
\[L(E;s) = \prod_{p \, \nmid \, \Delta(E_p)}\frac{1}{1-t_p p^{-s}+p^{1-2s}}\cdot \prod_{p \, \mid \, \Delta(E_p)}{}^{\prime}(1-p^{-s})\cdot \prod_{p \, \mid \, \Delta(E_p)}{}^{\prime \prime}(1+p^{-s})\]
gegeben, wobei wir mit $\prod{}^{\prime}$ als das Produkt \"uber allen Primzahlen $p$, wo $E$ split-multiplikative Reduktion hat und mit $\prod{}^{\prime \prime}$ als das Produkt \"uber allen Primzahlen $p$, wo $E$ nicht-split-multiplikative Reduktion hat, bezeichnen. Man definiert
\begin{equation*}
f_p = \begin{cases}0, \quad p \nmid \Delta(E_p), \\ 1, \quad p \mid \Delta(E_p), E_p^{sm}\otimes \mathbb{F}_p^{alg} \cong \mathbb{G}_m(\mathbb{F}_p^{alg}) \\ 2, \quad p \mid \Delta(E_p), E_p^{sm}\otimes \mathbb{F}_p^{alg} \cong \mathbb{G}_a(\mathbb{F}_p^{alg}), p \neq 2,3 \\ 2 + \delta_p, \quad p \mid \Delta(E_p), E_p^{sm}\otimes \mathbb{F}_p^{alg} \cong \mathbb{G}_a(\mathbb{F}_p^{alg}), p = 2,3 \end{cases}
\end{equation*}
wobei $\delta_p$ eine unver\"anderliche Gr\"o\ss e ist, die die wilde Verzweigung von der Wirkung der Tr\"agheitsgruppe $I_p \trianglelefteq \mathrm{Gal}(\mathbb{Q}^{alg}/\mathbb{Q})$ auf dem Tate-Modul $T_p(E)$ beschreibt (vgl. \cite[Appendix C. \S 16]{Sil09} f\"ur weitere Literaturen \"uber die technische Definition von $\delta_p$). Der Konduktor $N_{E/\mathbb{Q}}$ wird durch
\[N_{E/\mathbb{Q}} = \prod_{p \in \mathbb{P}}p^{f_p}.\]
definiert. Das unendliche Produkt $L(E;s)$ konvergiert f\"ur $\sigma > \frac{3}{2}$ wegen der Hasse-Schranke  $\mid t_p \mid < 2\sqrt{p}$. Daher hat man die Darstellung von $L(E;s)$ als eine Dirichlet-Reihe
\[L(E;s) = \sum_{n=1}^{+\infty}\frac{a_n}{n^s}, \quad \sigma > \frac{3}{2}.\]
Sei nun
\[\xi_E(s) = N_{E/\mathbb{Q}}^{s/2}(2\pi)^{-s}\Gamma(s)L(E;s).\]
Es ist bekannt, dass $L(E,s)$ sich analytisch auf die ganze Ebene $\mathbb{C}$ fortsetzen l\"asst und  $\xi_E(s)$ die Funktionalgleichung
\[\xi_E(s) = w \cdot \xi_E(2-s), \quad w = \pm 1\]
erf\"ullt. Denn diese Tatsache ist eine Folgerung aus der Taniyama-Shimura-Weil-Vermutung, die in einer ganzen Reihen von Arbeiten von Wiles, Taylor-Wiles, Breuil, Conrad, Diamond und Taylor bewiesen wurde. Wir werden hier nicht \"uber diese Vermutung diskutieren. Wir definieren nun die zu $E$ assoziierte globale erzeugende Funktion
\begin{equation}\label{e6}
B_{E/\mathbb{Q}}(s) \stackrel{defn}{=} \prod_{p \in \mathbb{P}, \, p \, \nmid \, \Delta(E_p)}\frac{(p^s-1)B_{E,p}(s)}{\# E(\mathbb{F}_p)} = \prod_{p \in \mathbb{P}, \, p \, \nmid \, \Delta(E_p)}\frac{1}{1-t_pp^{-s}+p^{1-2s}}\cdot \frac{1-p^{1-2s}}{1-p^{1-s}}.
\end{equation}
\begin{satz}\label{s10}
Sei $E/\mathbb{Q}$ eine elliptische Kurve. Dann hat man eine Funktionalgleichung
\[B_{E/\mathbb{Q}}(s) = \prod_{p \, \mid \, \Delta(E_p)}(1-p^{1-s})\cdot \prod_{p \, \mid \, \Delta(E_p)}(1-p^{1-2s})^{-1}\cdot A_{E/\mathbb{Q}}(s)\cdot \zeta(s-1)\cdot \zeta(2s-1)^{-1}, \quad \sigma > 2\]
wobei
\[A_{E/\mathbb{Q}}(s)= \frac{L(E;s)}{\prod_{p \, \mid \, \Delta(E_p)}{}^{\prime}(1-p^{-s})\cdot \prod_{p \, \mid \, \Delta(E_p)}{}^{\prime \prime}(1+p^{-s})}.\]
\end{satz}
\begin{proof}
Trivial.
\end{proof}
\begin{bem}
{
\rm
Die Funktion $A_{E/\mathbb{Q}}(s)$ ist die Globalisierung von der erzeugenden Funktion $A_{E}(p^{-s})$ der Folge $\{t_n\}_{n \in \mathbb{N}}$ in der rekursiven Formel \ref{e0}.
}
\end{bem}
\begin{satz}\label{s11}
Sei $E/\mathbb{Q}$ eine elliptische Kurve. Dann stellt die zu $E$ assoziierte globale erzeugende Funktion $B_{E/\mathbb{Q}}(s)$ auf $\{\mathfrak{R}(s) > 2\}$ eine analytische Funktion dar.
\end{satz}
\begin{proof}
Es gilt f\"ur $\sigma > 1$
\[\big \vert \frac{1}{\zeta(2s-1)} \big \vert = \big \vert \sum_{n=1}^{+\infty}\frac{\mu(n)}{n^{2s-1}} \big \vert \leq \sum_{n=1}^{+\infty}\frac{1}{n^{2\sigma - 1}} < \infty.\]
$\zeta(2s-1)$ hat keine Nullstelle auf $\sigma > 1$. Man hat
\[\# \{p \in \mathbb{P}: p \, \mid \, \Delta(E_p)\} < \infty\]
und $1-p^{1-2s} \neq 0$ f\"ur $\sigma \neq 1/2$. Man hat $1 \pm p^{-s} \neq 0$ f\"ur $\sigma \neq 0$. Das unendliche Produkt $L(E;s)$ ist analytisch im Bereich $\sigma > \frac{3}{2}$ und die Riemannsche Zetafunktion $\zeta(s-1)$ als unendliches Produkt  ist analytische auf $\sigma > 2$. Damit wurde der Satz bewiesen. 
\end{proof}
\begin{satz}\label{s12}
Sei $E/\mathbb{Q}$ eine elliptische Kurve. Die zu $E$ assoziierte globale erzeugende Funktion $B_{E/\mathbb{Q}}(s)$ ist analytisch fortsetzbar auf 
\[\mathbb{C}\setminus\{\mathfrak{R}(s) = 0\}\cup\{1/2 \leq \mathfrak{R}(s) < 1\}\cup \{ -n + \frac{1}{2}, n = 2,4,\cdots\} \cup \{2\}\]
und erf\"ullt dort auf dem Gebiet die Funktionalgleichung
\[B_{E/\mathbb{Q}}(s) = \prod_{p \, \mid \, \Delta(E_p)}(1-p^{1-s})\cdot \prod_{p \, \mid \, \Delta(E_p)}(1-p^{1-2s})^{-1}\cdot A_{E/\mathbb{Q}}(s)\cdot \zeta(s-1)\cdot \zeta(2s-1)^{-1}.\]
\end{satz}
\begin{proof}
Die Riemannsche Zetafunktion $\zeta(2s-1)$ hat triviale Nullstellen von der Form $-n+\frac{1}{2}$ f\"ur gerade positive nat\"urliche Zahl $n$ mit dem kritischen Streifen 
\[\{0 < \mathfrak{R}(2s-1) < 1\} = \{1/2 < \sigma < 1\}.\]
Die $L$-Funktion $L(E;s)$ ist fortsetzbar auf die ganze Ebene $\mathbb{C}$ und die Riemannsche Zetafunktion $\zeta(s-1)$ ist analytisch auf $\mathbb{C}\setminus \{2\}$. Die Faktoren $1-p^{1-2s}$ und $1 \pm p^{-s}$ sind nicht $0$ f\"ur $\sigma \neq 1/2$ und $\sigma \neq 0$. Daraus folgt die Behauptung.
\end{proof}
\begin{bem}
{\rm
Unter der Annahme der Riemannschen Vermutung ist $B_{E/\mathbb{Q}}(s)$ auf
\[\mathbb{C}\setminus\{\mathfrak{R}(s) = \frac{1}{2}\}\cup\{\mathfrak{R} = \frac{3}{4}\}\cup \{\mathfrak{R}(s) = 0\}\cup \{\mathfrak{R}(s) < 0, s \neq -n + \frac{1}{2}, n = 2,4,\cdots\} \cup \{2\}\]
fortsetzbar.
}
\end{bem}
\begin{hilfssatz}[\cite{Sch69}]\label{h4}
Es existieren Konstanten $\gamma_0 \geq 1$ und $\delta \in (0,1/2)$, so dass $\forall s = \sigma + it$ mit $\sigma \geq 1$, $\mid s-1 \mid \geq 1$ die Absch\"atzungen
\[\mid \zeta(\sigma + it) \mid \leq \gamma_0 \cdot (\mid t \mid + 2)^{\delta}, \quad \mid \zeta(\sigma + it) \mid^{-1} \leq \gamma_0 \cdot (\mid t \mid + 2)^{\delta}\]
gelten.
\end{hilfssatz} 
\begin{proof}
vgl. \cite[\S IV Satz 1.1, Satz 4.1, Satz 4.3]{Sch69}.
\end{proof}
\begin{satz}\label{s13}
Sei $E/\mathbb{Q}$ eine elliptische Kurve. F\"ur $\sigma > 2$ gibt es Konstanten $\gamma_1 \geq 1$ und $\delta \in (0,1/2)$, so dass die Absch\"atzung gilt
\[\mid B_{E/\mathbb{Q}}(s) \mid \leq \frac{36 \cdot \zeta(3) \cdot \zeta(3/2)}{\pi^2}\cdot \prod_{p \, \mid \, \Delta(E_p)}\big ( 1 + \frac{1}{p} \big ) \gamma_1^2(\mid t \mid + 2)^{2\delta}.\]
\end{satz}
\begin{proof}
F\"ur $\sigma > 3/2$ konvergiert das unendliche Produkt $L(E;s)$ und l\"asst sich dort als eine Dirichletsche Reihe
\[L(E;s) = \sum_{n=1}^{+\infty}\frac{a_n}{n^s}, \quad s = \sigma + it,\]
darstellen. Aus dem Satz von Hasse hat man
\[\big \vert \frac{a_n}{n^s} \big \vert = \mathcal{O}(\frac{1}{n^{\sigma - 1/2}}).\]
Daher hat man f\"ur $\sigma > 2$
\[\mid L(E;s) \mid = \sum_{n=1}^{+\infty}\mathcal{O}(\frac{1}{n^{\sigma-1/2}}) \leq \sum_{n=1}^{+\infty}\frac{1}{n^{3/2}} = \zeta(3/2).\]
Es gilt f\"ur $\sigma > 2$
\[\prod_{p \, \mid \, \Delta(E_p)}\mid 1 - p^{1-s} \mid \leq \prod_{p \, \mid \, \Delta(E_p)}  1 + p^{1-\sigma} \stackrel{\sigma > 2}{\leq} \prod_{p \, \mid \, \Delta(E_p)}(1+\frac{1}{p}).\]
Da es gilt f\"ur $\sigma > 2$
\[\mid 1 - \frac{1}{p^{2s-1}} \mid \geq 1 - \frac{1}{p^{2\sigma-1}} \geq 1 - \frac{1}{p^3},\]
hat man 
\[\big \vert \prod_{p \, \mid \, \Delta(E_p)} \frac{1}{1-p^{1-2s}} \big \vert \leq \prod_{p \, \mid \, \Delta(E_p)}\frac{1}{1-\frac{1}{p^3}} \leq \prod_{p \in \mathbb{P}}\frac{1}{1-\frac{1}{p^3}} = \zeta(3).\]
Nun sch\"atzen wir ab 
\[\big\vert \prod_{p \, \mid \, \Delta(E_p)}{}^{\prime}(1-p^{-s}) \big\vert = \prod_{p \, \mid \, \Delta(E_p)}{}^{\prime}\big\vert 1 - p^{-s}\big\vert \geq \prod_{p \,\mid \,\Delta(E_p)}{}^{\prime}1-\frac{1}{p^{\sigma}}.\]
Daraus folgt
\[\big\vert \prod_{p \, \mid \, \Delta(E_p)}{}^{\prime}1-p^{-s}\big\vert^{-1} \leq \prod_{p \, \mid \, \Delta(E_p)}{}^{\prime}(1-\frac{1}{p^{\sigma}})^{-1} \leq \zeta(\sigma) \leq \zeta(2) = \frac{\pi^2}{6}.\]
Das Gleiche gilt auch f\"ur $\mid \prod_{p \, \mid \, \Delta(E_p)} 1+p^{-s} \mid \leq \frac{\pi^{2}}{6}$. Aus dem Hilfssatz \ref{h4} gibt es Konstanten $\gamma_0, \gamma_0^{\prime} \geq 1$ und $\delta \in (0,1/2)$, so dass es f\"ur $\sigma > 2$ gilt
\[\mid \zeta(s-1) \mid \leq \gamma_0 \cdot (\mid t \mid + 2)^{\delta}, \quad \mid \zeta(2s-1) \mid^{-1} \leq \gamma_0^{\prime} \cdot (\mid t \mid + 2)^{\delta}.\]
Man nimmt dann $\gamma_1 = \max\{\gamma_0,\gamma_0^{\prime}\}$. Damit wurde der Satz bewiesen.
\end{proof}
\begin{satz}\label{s14}
Seien $E/\mathbb{Q}$ eine elliptische Kurve und $B_{E/\mathbb{Q}}(s)$ die dazugeh\"orige globale erzeugende Funktion. Als eine analytische Funktion auf dem Gebiet
\[\mathbb{C}\setminus\{\mathfrak{R}(s) = 0\}\cup\{1/2 \leq \mathfrak{R}(s) < 1\}\cup \{-n + \frac{1}{2}, n = 2,4,\cdots\} \cup \{2\}\]
hat $B_{E/\mathbb{Q}}(s)$ eine Nullstelle an dem Punkt $s=1$. D.h
\[\lim_{s\to 1+}B_{E/\mathbb{Q}}(s) = 0.\]
\end{satz}
\begin{proof}
Die Behauptung folgt aus der bekannten Tatsache, dass $\zeta(2s-1)$ hat einen einfachen Pol an dem Punkt $s=1$ hat, d.h
\[\lim_{s\to 1+}\frac{1}{s-1}\zeta(2s-1) = 1.\]
\end{proof}
\begin{bem}
{\rm
Es wurde von Birch und Swinnerton-Dyer vermutet, dass der Zentralwert $L(E;1)$ der $L$-Funktion gleich $0$ ist. Etwas pr\"aziser wurde es vermutet, dass $L(E;s)$ eine Nulstelle der Vielfachheit $\mathrm{rang}(E(\mathbb{Q}))$ an der Stelle $s=1$ besitzt. Wenn man die Richtigkeit der Birch-Swinnerton-Dyer-Vermutung annimmt, dann hat man $\mathrm{ord}_{s\to 1+}B_{E/\mathbb{Q}}(s) = \mathrm{rang}(E(\mathbb{Q}))+1$.
}
\end{bem}
\begin{bem}
{\rm
Die folgende Anmerkung hat keine mathematische, sondern nur eine heuristische Bedeutung. Man kann versuchen, $s=1$ in die Darstellung von $B_{E/\mathbb{Q}}(s)$ als ein unendliches Produkt in der Gleichung \ref{e6} einzusetzen, obwohl das Produkt dort in $s=1$ auf gar keinen Fall konvergiert. Dank dem Satz \ref{s14} wissen wir, dass die analytische Forsetzung $B_{E/\mathbb{Q}}(1) = 0$ erf\"ullt. So erh\"alt man $\prod_{p \, \nmid \, \Delta(E_p)}\frac{\# E(\mathbb{F}_p)}{(p-1)} \to +\infty$. Diese Beobachtung ist \"ahnlich wie die heuristische Beobachtung mit dem Zentralwert $L(E;1)$ der $L$-Funktion $L(E;s)$.  
}
\end{bem}
\begin{satz}\label{s15}
Seien $E/\mathbb{Q}$ eine elliptische Kurve und $B_{E/\mathbb{Q}}(s)$ die zu $E$ assoziierte globale erzeugende Funktion. Man hat eine Funktionalgleichung
\[B_{E/\mathbb{Q}}(1-s) = B_{E/\mathbb{Q}}(s)\cdot \zeta(2s-1)^{-1}\cdot \zeta(s) \cdot \frac{\prod_{p \, \mid \, \Delta(E_p)}1-\frac{1}{p^s}}{\prod_{p \, \mid \, \Delta(E_p)}1-\frac{1}{p^{2s-1}}}.\]
auf dem Gebiet
\[\mathbb{C}\setminus\{\mathfrak{R}(s) = 0\}\cup\{1/2 \leq \mathfrak{R}(s) < 1\}\cup \{-n + \frac{1}{2}, n = 2,4,\cdots\} \cup \{2\}\]
\end{satz}
\begin{proof}
Die Behauptung folgt aus der Funktionalgleichung \ref{e5} des lokalen Faktors $B_{E,p}(s)$. 
\end{proof}
\noindent Die Dirichletsche Reihe von $L(E;s)$ ist gegeben durch
\[L(E;s) = \sum_{n=1}^{+\infty}\frac{a_n}{n^{s}}, \quad s = \sigma + it, \sigma > \frac{3}{2}\]
mit
\[a_p = \begin{cases}t_p, \quad p \nmid \Delta(E_p) \\ 1, \quad p \mid \Delta(E_p), \, \text{split-multiplikativ} \\ -1, \quad p \mid \Delta(E_p), \, \text{nicht-split-multiplikativ} \\ 0, \quad p \mid \Delta(E_p), \, \text{additiv}\end{cases}\]
F\"ur $k \geq 1$ definiert man die Fourier-Koeffizienten rekursiv
\[\begin{cases} a_{p^{k+1}} = a_p \cdot a_{p^k} - p\cdot a_{p^{k-1}}, \quad p \nmid \Delta(E_p) \\ a_{p^k} = a_p^k, \quad p \mid \Delta(E_p)\end{cases}\]
und $a_{n\cdot m} = a_n \cdot a_m$, falls $\mathrm{ggT}(n,m) = 1$. Seien $\sigma(n)$ und $d(n)$ die zahlentheoretischen Funktionen
\[\sigma(n) = \sum_{d \, \mid \, n}d, \quad d(n) = \sum_{d \, \mid \, n}1.\]
F\"ur zwei zahlentheoretische Funktionen $f$ und $g$ bezeichnen wir $f*g$ die Dirichletsche Faltung.
\begin{satz}\label{s16}
Seien $E/\mathbb{Q}$ eine elliptische Kurve und $B_{E/\mathbb{Q}}(s)$ die zu $E$ assoziierte globale erzeugende Funktion. F\"ur $\mathfrak{R}(s) > 2$ nimmt $B_{E/\mathbb{Q}}(s)$ die Form vom Produkt einer Dirichletschen Reihe mit einem endlichen Produkt von lokalen Faktoren
\[B_{E/\mathbb{Q}}(s) = \frac{\prod_{p \, \mid \, \Delta(E_p)}(1-p^{1-s})\cdot \prod_{p \, \mid \, \Delta(E_p)}(1-p^{1-2s})^{-1}}{\prod_{p \, \mid \, \Delta(E_p)}{}^{\prime}(1-p^{-s})\cdot \prod_{p \, \mid \, \Delta(E_p)}{}^{\prime \prime}(1+p^{-s})} \sum_{n=1}^{+\infty}\frac{b_n}{n^s},\]
deren Fourier-Koeffizienten $b_n$'s durch
\[b_n = \sum_{k_5\cdot k_6 = n}a_{k_6}(\sum_{k_3\cdot k_4 = k_5}\sigma(k_4)\cdot d(k_4)\cdot (\sum_{k_1\cdot k_2 = k_3}\mu * \mu(k_1) \cdot k_2 \cdot \mu(k_2))).\]
dargestellt werden k\"onnen.
\end{satz}
\begin{proof}
F\"ur $\mathfrak{R}(s) > 2$ hat man die Ramanujan-Formel (vgl. \cite[Chap. I 1.3]{Titch51})
\[\frac{\zeta(s)^2 \zeta(s-1)^2}{\zeta(2s-1)} = \sum_{n=1}^{+\infty}\frac{\sigma(n)\cdot d(n)}{n^s}.\]
Man hat f\"ur $\mathfrak{R}(s) > 1$
\[\frac{1}{\zeta(s)^2} = \sum_{n=1}^{+\infty}\frac{\mu(n)}{n^s}\sum_{n=1}^{+\infty}\frac{\mu(n)}{n^s} = \sum_{n=1}^{+\infty}\frac{\mu*\mu(n)}{n^s}.\]
Damit ist
\[\frac{1}{\zeta(s)^2\cdot \zeta(s-1)} = \sum_{n=1}^{+\infty}\frac{\mu*\mu(n)}{n^s}\cdot \sum_{n=1}^{+\infty}\frac{n \cdot \mu(n)}{n^s} = \sum_{n=1}^{+\infty}\frac{\sum_{k \, \mid \, n}\mu * \mu(k) \cdot \frac{n}{k} \cdot \mu(\frac{n}{k})}{n^s}.\]
Die Behauptung folgt nun unmittelbar aus dem Satz \ref{s10}.
\end{proof}
\begin{satz}\label{s17}
Sei $E/\mathbb{Q}$ eine elliptische Kurve. $B_E(s)$ hat einen einfachen Pol an der Stelle $s=2$ mit
\[\mathrm{Res}_{s=2}(B_{E/\mathbb{Q}}(s)) = \frac{\prod_{p\,\mid\,\Delta(E_p)}(1-\frac{1}{p})(1-\frac{1}{p^3})^{-1}}{\prod_{p \, \mid \, \Delta(E_p)}^{}{\prime}(1-\frac{1}{p^2})\cdot \prod_{p\,\mid\,\Delta(E_p)}{}^{\prime \prime}(1+\frac{1}{p^2})}\cdot \frac{L(E,2)}{\zeta(3)}. \]
\end{satz}
\begin{proof}
Es gilt $\mathrm{Res}_{s=2}\zeta(s-1) = 1$ und $\mathrm{ord}_{s=2}\zeta(s-1) = 1$. Daraus folgt die Behauptung.
\end{proof}
\begin{folgerung}\label{f2}
Sei $E/\mathbb{Q}$ eine elliptische Kurve. Dann hat man eine asymptotische Formel
\[B_{E/\mathbb{Q}}(s) \sim \frac{\prod_{p\,\mid\,\Delta(E_p)}(1-\frac{1}{p})(1-\frac{1}{p^3})^{-1}}{\prod_{p \, \mid \, \Delta(E_p)}^{}{\prime}(1-\frac{1}{p^2})\cdot \prod_{p\,\mid\,\Delta(E_p)}{}^{\prime \prime}(1+\frac{1}{p^2})}\cdot \frac{L(E,2)}{\zeta(3)} \cdot \frac{1}{s-2}, \quad s \to 2+ \]
\end{folgerung}
\begin{proof}
Der Beweis folgt aus dem Satz \ref{s17}.
\end{proof}
\begin{bem}
{\rm
Die Fourier-Koeffizienten $a_n$ in der Dirichletschen Reihe $L(E;s)$ und $b_n$ in der globalen erzeugenden Funktion $B_{E/\mathbb{Q}}(s)$ haben i. Allg. abwechselndes Vorzeichen. Daher d\"urfen wir nicht die Taubers\"atze anwenden, um das asymptotische Verhalten von $\lim_{x \to +\infty}\frac{1}{x}\cdot \sum_{n \leq x}b_n$ zu erschlie\ss en. Im Fall der M\"obius-Funktion kann man immer $1$ dazu addieren und die Formel
\[\sum_n (1+\mu(n))n^{-1-\sigma} = \zeta(1+\sigma) + \zeta(1+\sigma)^{-1}\]
anwenden. Leider k\"onnen wir diesen Trick nicht f\"ur $\{b_n\}$ i. Allg. verwenden. Wegen der Hasse-Schranke hat man nur
\[\mid a_{n} \mid \leq \sqrt{n}\cdot d(n).\]
}
\end{bem}
\noindent Riemann hat folgende Formel in dem Beweis seiner Funktionalgleichung erhalten
\[\Gamma(\frac{s}{2})\pi^{-\frac{s}{2}}\zeta(s) = \frac{1}{s(s-1)} + \int_1^{+\infty}\psi(x)(x^{s/2-1}+x^{-1/2-s/2})dx,\]
wobei
\[\psi(x) = \sum_{n=1}^{+\infty}e^{-n^2\pi x}, \quad x > 0.\]
Die Riemannsche $\xi$-Funktion wird definiert durch
\[\xi(s) = \frac{1}{2}s(s-1)\pi^{-s/2}\Gamma(\frac{s}{2})\zeta(s).\]
Man setzt $H(z) = \xi(\frac{1}{2}+iz)$ f\"ur $z \in \mathbb{C}$. Bekanntlich sind $\xi(s)$ und $H(z)$ ganze Funktionen der Ordnung $1$ (vgl. \cite[Thm. 2.12]{Titch51}). Unter der Ordnung einer ganzen Funktion $f: \mathbb{C} \rightarrow \mathbb{C}$ versteht man die Gr\"o\ss e
\[\rho = \limsup_{r \to +\infty}\frac{\log(\log \norm{f}_{\infty,B_r})}{\log r}, \quad \norm{f}_{\infty,B_r} = \sup_{z \in B_r}\mid f(z) \mid,\]
wobei $B_r$ der Ball vom Radius $r$ ist. $H(z)$ ist reell f\"ur $z\in \mathbb{R}$. Aus der Funktionalgleichung von Riemann hat man
\[H(z) = \frac{1}{2} - (z^2+\frac{1}{4})\int_1^{+\infty}\psi(x)x^{-3/4}\cos(\frac{1}{2}z\log x)dx.\]
Aus der Theorie der Theta-Reihen, auch bekannt unter dem Name Poisson-Summenformel, hat man
\[\sqrt{x}\cdot(2\psi(x)+1) = 2\cdot \psi(\frac{1}{x})+1.\]
Man setzt
\[\varphi(x) = x^{1/4}(2x\frac{d^2}{dx^2}\psi(x)+3\frac{d}{dx}\psi(x)), \quad \Phi(t) = 2e^{2t}\varphi(e^{2t}), \quad t \in \mathbb{R}.\]
So hat man
\[H(z) = \int_{\mathbb{R}}\Phi(u)e^{izu}du\]
mit
\[\Phi(u) = 2\cdot \sum_{n=1}^{+\infty}(2n^4\pi^2e^{9u/2}-3n^2\pi e^{5u/2})e^{-n^2\pi e^{2u}}.\]
Die Nullstellen von $\zeta(s)$ auf der kritischen Gerade $\sigma = \frac{1}{2}$ sind genau die reellen Nullstellen von $H(z)$. Hardy hat bewiesen, dass $H(z)$ unendlich viele $\mathbb{R}$-Nullstellen hat (vgl. \cite[Chap X. 10.2]{Titch51}). Folglich haben wir
\begin{satz}\label{s18}
Seien $E/\mathbb{Q}$ eine elliptische Kurve und $B_{E/\mathbb{Q}}(s)$ die zu $E$ assoziierte globale erzeugende Funktion. Als eine analytische Funktion auf dem Gebiet
\[\mathbb{C}\setminus\{\mathfrak{R}(s) = 0\}\cup\{1/2 \leq \mathfrak{R}(s) < 1\}\cup \{-n + \frac{1}{2}, n = 2,4,\cdots\} \cup \{2\}\]
hat $B_{E/\mathbb{Q}}(s)$ unendlich viele Nullstellen auf der Gerade $\sigma = -\frac{1}{2}$.
\end{satz}
\noindent Im Zusammenhang mit Kongruenzformel der rationalen Punkte hat man folgendes Ergebnis von Deuring (vgl. \cite{Deu53}).
\begin{satz}[Deuring]\label{s19}
Seien $E/\mathbb{Q}$ eine elliptische Kurve mit komplexer Multiplikation und $x \geq 2$ eine Zahl. Dann hat man eine asymptotische Formel
\[\sum_{\mathclap{\substack{p \leq x,\, p \nmid \Delta(E_p)\\ \#E(\mathbb{F}_p) \equiv 1 \mod p}}}1 \sim \frac{1}{2}\int_2^x \frac{du}{\log u} + \underline{\mathcal{O}}_{x \to +\infty}(x\cdot \exp(-\frac{1}{200}\sqrt{\log(x)}), \quad x \to +\infty.\]
\end{satz}
\begin{bem}
Die Lang-Trotter-Vermutung besagt, dass es eine Konstante $C=C(E)>0$ gibt, so dass folgende asymptotische Formel gilt  (vgl. \cite[p. 144]{Sil09})
\[\#\{p: p \nmid \Delta(E_p), p \leq x, E_p \quad \text{supersingul\"ar} \} \sim C\cdot\frac{\sqrt{x}}{\log x}, \quad x \to +\infty,\]
falls $E$ keine komplexe Multiplikation besitzt. Wir erw\"ahnen an dieser Stelle noch eine Vermutung von Koblitz in \cite[Conj. A]{Kob88}. F\"ur eine elliptischen Kurve $E/\mathbb{Q}$ mit $\Delta_E \in \mathbb{Z}$, die nicht $\mathbb{Q}$-isogen zu einer Kurve mit nicht-trivialer Torsiongruppe $E_{\mathrm{tor}}$ der $\mathbb{Q}$-Punkte ist, hat Koblitz vermutet
\[\sum_{\mathclap{\substack{p \leq x,\, p \nmid \Delta(E_p)\\ \#E(\mathbb{F}_p) \in \mathbb{P}}}}1 \sim C\cdot \frac{x}{\log^2 x}, \quad x \to +\infty,\]
falls $E$ keine komplexe Multiplikation besitzt, wobei $C$ eine positive Konstante ist. Es w\"are interessant, wenn man einen Zusammenhang zwischen den Funktionen 
\[\sum_{\mathclap{\substack{p \leq x,\, p \nmid \Delta(E_p)\\ \#E(\mathbb{F}_p) \equiv 1 \mod p}}}1\] 
und $B_{E/\mathbb{Q}}(s)$ erkennen k\"onnte. Denn $B_{E/\mathbb{Q}}(s)$ verf\"ugt eine Mischung aus der Riemannschen Zetafunktion, die die Verteilung von Primzahlen liefert, sowie aus der $L$-Funktion der elliptischen Kurven, die die arithmetischen Eigenschaften von der Kurve nachweist.
\end{bem}
\noindent Im folgenden wollen wir Modulformen betrachten. Viele S\"atze lassen sich leicht verallgemeinern aber wir beschr\"anken uns nur auf Spitzenformen und zwar auf der ganzen Modulgruppe $\mathbf{SL}_2(\mathbb{Z})$. Man kann nat\"urlich auf Level $N$ ohne gro\ss e Anstrengung verallgemeinern. Wir ersparen uns an dieser Stelle die Einzelheiten. Die Hauptmotivation f\"ur uns ist die klassische Arbeit von G. P\'olya \cite{Polya26}. Sei nun $f \in S_k(\Gamma_1)$ eine Spitzenform vom Gewicht $k$ mit $\Gamma_1 = \mathbf{SL}_2(\mathbb{Z})$. Als eine Fourier-Reihe hat man 
\[f = \sum_{n=1}^{+\infty}a_n \exp(2\pi i n z), \quad a_n = \mathcal{O}(n^{k/2}), z \in \mathbb{H} = \{z \in \mathbb{C}: \mathfrak{Im}(z) > 0\}.\]
Die $L$-Funktion von Hecke 
\[L(f;s) = \sum_{n=1}^{+\infty}\frac{a_n}{n^s}\]
konvergiert absolut und lokal gleichm\"a\ss ig auf dem Gebiet $\mathfrak{R}(s) > \frac{k}{2}+1$. Ferner, $L(f;s)$ ist fortsetzbar auf die ganze Ebene $\mathbb{C}$ und erf\"ullt die Funktionalgleichung
\[(2\pi)^{-s}\Gamma(s)L(f;s) = (-1)^k(2\pi)^{s-k}\Gamma(k-s)L(f;k-s).\]
Man hat die Mellin-Transform
\[\xi_f(s) \stackrel{defn}{=} (2\pi)^{-s}\Gamma(s)L(f;s) = \int_1^{+\infty}\phi(t)(t^s + (-1)^{k/2} t^{k-s})\frac{dt}{t},\]
wobei
\[\phi(t) = \sum_{n=1}^{+\infty}a_n \exp(-2\pi n t), \quad t > 0.\]
Es gilt
\[\phi(\frac{1}{t}) = (-1)^{k/2}t^{k}\phi(t).\]
Sei nun $f \in S_2$ eine Spitzenform vom Gewicht $2$. Nach zweimal partiellen Integrationen mit den Nebenbedingungen $\phi(1) = \phi(1)^{\prime} = 0$, die man aus der Funktionalgleichung f\"ur $\phi(t)$ oben leicht ermitteln kann,  erh\"alt man eine Funktion
\[\xi_f(s) = 2\int_1^{+\infty}\phi(t)t\sinh((s-1)\log t)dt.\]
Man setzt $t=e^u$ und definiert $H_f(z) = \xi_f(1+z)$. So erh\"alt man die Integraldarstellung
\[H_f(z) = \int_0^{+\infty}\Phi_f(u)\sinh(zu)du\]
mit 
\[\Phi_f(u) = 2e^{2u}\sum_{n=1}^{+\infty}a_n e^{-2\pi n e^u}.\]
Es ist offensichtlich, dass $H_f(z)$ reell f\"ur $z \in \mathbb{R}$ ist. Ferner, $H_f(z)$ hat die Potenzreihe-Entwicklung
\[H_f(z) = \sum_{n=0}^{+\infty}\frac{z^{2n+1}}{(2n+1)!}\hat{b}_{2n+1}, \quad \hat{b}_{2n+1} = \int_{0}^{+\infty}\Phi_f(u)u^{2n+1}du.\]
Wegen der Funktionalgleichung von $\xi_f(s)$ ist $H_f(z)$ eine ungerade Funktion
\[H_f(z) = -H_f(-z).\]
F\"ur $\lambda \in \mathbb{R}$ definiert man
\[H_{f,\lambda}(z) = \int_0^{+\infty}e^{\lambda u}\Phi_f(u)\sinh(zu)du.\]
Diese Funktion erf\"ullt die Wellengleichung von D'Alembert
\[\partial_{\lambda \lambda}H_{f,\lambda}(z) = \partial_{zz}H_{f,\lambda}(z).\]
\begin{satz}\label{s20}
Sei $f \in S_2$ eine Spitzenform vom Gewicht $2$. Dann ist $H_f(z)$ eine ganze Funktion der Ordnung $1$.
\end{satz}
\begin{proof}
Wegen der Funktionalgleichung $\xi_f(s) = -\xi_f(2-s)$ reicht es, $\mid \xi_f(s) \mid$ auf dem Gebiet $\sigma > 2$ abzusch\"atzen. F\"ur eine Spitzenform haben wir auf dem Gebiet $\sigma > 2$ wegen dem Hilfssatz \ref{h1}
\begin{equation*}
\begin{split}
\sum_{n=1}^{+\infty}a_nn^{-s} = \lim_{x \to +\infty} \sum_{\log n \leq x}a_n \cdot \exp(-s\log n) \\ 
= \lim_{x \to +\infty}((\sum_{\log n \leq x}a_n)\cdot e^{-sx} + s\cdot \int_0^x \sum_{\log n \leq u} a_n \cdot e^{-su}du) = s\cdot \int_0^{+\infty}(\sum_{n \leq e^u}a_n)e^{-su}du.
\end{split}
\end{equation*}
Da $a_n = \mathcal{O}(n)$ ist, hat man
\[\mid s\cdot \int_0^{+\infty}(\sum_{n \leq e^u}a_n)e^{-su}du \mid \leq \mid s \mid \int_0^{+\infty}e^{-(\sigma-2)u}du = \frac{1}{\sigma-2}\mid s \mid.\]
Somit konvergiert das Integral $\int_0^{+\infty}(\sum_{n \leq e^u}a_n)e^{-su}du$ auf dem Gebiet $\sigma > 2$ und $L(f;s)$ hat dort die Asymptotik
\[L(f;s) = \mathcal{O}(\mid s \mid), \quad \mid s \mid \to +\infty.\]
Nun gilt es $(2\pi)^{-s} = \mathcal{O}(\exp(-\sigma\log(2\pi)))$. F\"ur die Gammafunktion gilt die folgende Absch\"atzung bekanntlich
\[\mid \Gamma(s) \mid \leq \exp(\mathrm{const} \cdot \mid s\mid \cdot \log \mid s\mid).\]
D.h. $H_f(z)$ hat h\"ochstens die Ordnung $1$. Aber wegen der Stirlingschen Formel
\[\log \Gamma(s) = (s-\frac{1}{2})\log(s) - s + \frac{1}{2}\log(2\pi) - \int_0^{+\infty}\frac{t-[t]-1/2}{z+t}dt\]
hat man die asymptotische Formel
\[\Gamma(s) \sim \sqrt{2\pi}s^{s-1/2}e^{-s}, \mid s \mid \to +\infty, \, s \in \mathbb{C}_{-},\]
wobei wir den Hauptzweig f\"ur $\log s$ in der geschlitzten Ebene
\[\mathbb{C}_{-} = \mathbb{C} - \{a \in \mathbb{R}: a \leq 0\}\]
nehmen. Daraus folgt
\[\log \Gamma(s) \sim s\log(s), \quad \mid s \mid \to +\infty, \, s \in \mathbb{C}_{-}.\]
So muss die ganze Funktion $H_f(z)$ die Ordnung $1$ haben. Damit wurde der Satz bewiesen.
\end{proof}
\begin{satz}\label{s21}
Sei $f \in S_2$ eine Spitzenform vom Gewicht $2$. Dann hat die Funktion $H_f(z)$ unendlich viele Nullstellen .
\end{satz}
\begin{proof}
Aus der Theorie von Faktorisierung der ganzen Funktionen von Hadamard und Weierstra\ss \, hat man die Darstellung f\"ur $z,A,B,\gamma_n \in \mathbb{C}$
\[H_f(z) = (z+1)^r \cdot e^{A+Bz}\prod_{n=1}^{+\infty}(1-\frac{z+1}{\gamma_n})e^{(1+z)/\gamma_n},\]
weil man die Absch\"atzung $\mid \xi_f(1+z) \mid \leq e^{c\mid 1 + z\mid}$ f\"ur ein festes $c>0$ hinreichend gro\ss \, hat (siehe den Satz \ref{s23} unten f\"ur eine asymptotische Absch\"atzung von einer oberen Schranke). Es ist bekannt, dass die Reihe $\sum_{n=1}^{+\infty}\mid \gamma_n \mid^{-1}$ divergiert und die Reihe $\sum_{n=1}^{+\infty}\mid \gamma_n \mid^{-1-\varepsilon}$ f\"ur jedes $\varepsilon > 0$ konvergiert. So muss das Produkt in der Darstellung von $H_f(z)$ oben unendlich viele Faktoren in dem Produkt
\[\prod_{n=1}^{+\infty}(1-\frac{1+z}{\gamma_n})e^{(1+z)/\gamma_n}\]
erhalten. Daraus folgt die Behauptung. 
\end{proof}
\begin{hilfssatz}[Lemma von Riemann-Lebesgue]\label{h5}
Es sei $f\in L^1(0,+\infty)$ integrierbar. F\"ur die Laplace-Transform von $f$ gilt
\[\int_0^{+\infty}f(t)e^{-tz}dt \to 0, \quad \mid z \mid \to +\infty,\, \mathfrak{R}(z) \geq 0.\]
\end{hilfssatz}
\begin{satz}\label{s22}
Sei $f \in S_2$ eine Spitzenform vom Gewicht $2$ mit 
\[H_f(z) = \int_0^{+\infty}\Phi_f(u)\sinh(zu)du.\]
Dann ist $\Phi_f(u)$ eine gerade Funktion.
\end{satz}
\begin{proof}
Der Beweis verl\"auft bekanntlich wie aus der klassischen Arbeit von G. P\'olya \cite{Polya26}. Angenommen, $\Phi_f(u)$ w\"are keine gerade Funktion. Dann g\"abe es eine nat\"urliche Zahl $n$, so dass
\[\Phi_f{}^{\prime}(0) = \Phi_f{}^{\prime \prime \prime}(0) = \cdots = \Phi_f^{(2n-3)}(0) = 0,\]
aber 
\[\Phi_f^{(2n-1)}(0) \neq 0.\] 
Wir betrachten die Funktion $H_f(y+ix)$ mit $x,y \in \mathbb{R}$. Durch wiederholte partielle Integration hat man f\"ur $x,y \in \mathbb{R}$
\[H_f(y+ix) = (-1)^n\frac{\Phi_f^{(2n-1)}(0)}{(y+ix)^{2n}} + \frac{(-1)^{n+1}}{(y+ix)^{2n+1}}\int_0^{+\infty}\Phi_f^{(2n+1)}(u)\cosh((y+ix)u)du.\]
Die Funktion
\[\Phi_f(u) = 2e^{2u}\sum_{n=1}^{+\infty}a_n e^{-2\pi n e^u}\]
ist eine Schwartz-Funktion, d.h. sie f\"allt drastisch schnell. Zun\"achst haben wir
\begin{equation*}
\begin{split}
\lim_{x \to +\infty}(y+ix)^{2n}H_f(y+ix) = (-1)^n\Phi_f^{(2n-1)}(0) + \\ \lim_{x \to +\infty}\frac{(-1)^{n+1}}{y+ix}\int_0^{+\infty}\Phi_f^{(2n+1)}(u)\cosh((y+ix)u)du.
\end{split}
\end{equation*}
Es gilt $0 \neq (-1)^n\Phi_f^{(2n-1)}(0) < \infty$. Nach dem Lemma von Riemann-Lebesgue (vgl. Hilfssatz \ref{h5}) haben wir
\[\int_0^{+\infty}\Phi_f^{(2n+1)}(u)e^{\pm (y+ix)u}du \to 0, \quad x \to \pm \infty.\]
Daraus folgt
\[\lim_{x \to +\infty}\frac{1}{y+ix}\int_0^{+\infty}\Phi_f^{(2n+1)}(u)\cosh((y+ix)u)du = 0.\]
Dies bedeutet, dass es
\[0 \neq \lim_{x \to +\infty}(y+ix)^{2n}H_f(y+ix)\]
gelten muss. Folglich ist $H_f(y+ix) \neq 0$ f\"ur $x$ hinreichend gro\ss, was aber zum Widerspruch zu dem Satz \ref{s21} f\"uhrt. So muss $\Phi_f(u)$ eine gerade Funktion sein. Damit wurde der Satz bewiesen.
\end{proof}
\begin{bem}
{\rm
Die Anwendung von dem Lemma von Riemann-Lebesgue in dem Beweis vom Satz \ref{s22} oben ist entbehrlich. Es gen\"ugt, zu beobachten, dass $\Phi_f^{(2n+1)}(u)$ eine Schwartz-Funktion ist und man hat daher eine Absch\"atzung
\[\big \vert \int_0^{+\infty}\Phi_f^{(2n+1)}(u)\cosh((y+ix)u)du \big \vert \, \leq \int_0^{+\infty}\mid \Phi_f^{(2n+1)}(u)\mid \cosh(yu)du,\]
was uns ohne Weiteres den gew\"unschten Grenzwert
\[\lim_{x \to +\infty}\frac{1}{y+ix}\int_0^{+\infty}\Phi_f^{(2n+1)}(u)\cosh((y+ix)u)du = 0\]
liefert.
}
\end{bem}
\noindent F\"ur eine Funktion $F(z)$ bezeichnen wir mit 
\[N(F(z)) = \{z \in \mathbb{C}: F(z) = 0\}\] 
die Menge ihrer Nullstellen. Nun definieren wir die Konstante von De Bruijn-Newman f\"ur eine Spitzenform $f$ vom Gewicht $2$
\[\Lambda_f = \inf_{\lambda \in \mathbb{R}}\{\lambda: N(H_{f,\lambda}(z)) \subset i\cdot\mathbb{R}\}.\] 
Die Menge der nicht-trivialen Nullstellen der Hecke-$L$-Funktion $L(f;s)$ einer gemeinsamen Eigenform $f \in S_k(\Gamma_0(N))^{new}$ des Hecke-Operators vom Gewicht $k$ und Level $N$ liegt in dem kritischen Streifen
\[\frac{k-1}{2} \leq \sigma \leq \frac{k+1}{2}.\]
Im Fall vom Gewicht $k=2$ ist die kritische Gerade $\sigma = 1$. Falls die Konstante von De Bruijn-Newman die Ungleichung $\Lambda_f \leq 0$ erf\"ullt, dann hat $H_f(z)$ nur $i\cdot\mathbb{R}$-Nullstellen. Dies ist \"aquivalent zur Vermutung, dass die nicht-triviale Nullstellen von $L(f;s)$ sich auf der kritischen Gerade $\sigma = 1$ konzentrieren. Die Funktion $H_f(z)$ hat eine triviale Nullstelle $z=0$. Man beachtet, dass $H_f(z)$ sich auch als
\[H_f(z) = \int_0^{+\infty}zu\cdot \Phi_f(u) \cdot \prod_{n=1}^{+\infty}(1+\frac{z^2u^2}{n^2\pi^2})du\]
schreiben l\"asst. Da es
\[H_f(\bar{z}) = \int_0^{+\infty}\Phi_f(u)\sinh(\bar{z}u)du = \int_0^{+\infty}\Phi_f(u)\overline{\sinh(zu)}du = \overline{H_f(z)}\]
und $H_f(z) = -H_f(-z)$ gelten, k\"onnen wir folglich schlie\ss en, dass $-z$ und $\bar{z}$ auch Nullstellen von $H_f(z)$, falls $z$ eine Nullstelle von $H_f(z)$ ist.
\begin{hilfssatz}\label{h6}
Sei $g: [a,x] \rightarrow \mathbb{R}$ eine monotone nicht-wachsende Funktion. Sind $a$ und $b$ ganz, so ist
\[\sum_{a < n \leq b}g(n) \leq \int_a^b g(u)du.\]
\end{hilfssatz}
\begin{proof}
Wegen der Monotonie von $g$ hat man die Absch\"atzung
\[\int_a^bg(u)du = \sum_{m=a+1}^b \int_{m-1}^m g(u)du \geq \sum_{m=a+1}^b g(m)\int_{m-1}^m du = \sum_{a < n \leq b}g(n) \]
\end{proof}
\begin{satz}\label{s23}
Sei $f \in S_2$ eine Spitzenform vom Gewicht $2$. Dann hat man eine asymptotische Formel
\[H_f(z) = \mathcal{O}(\int_0^{+\infty}e^u\cosh(xu)du), \quad z = x + iy, \mid z \mid \to +\infty.\]
\end{satz}
\begin{proof}
Wir haben
\[\mid H_f(z) \mid \, = \, \mid \int_0^{+\infty}\Phi_f(u)\sinh(zu)du \mid \, \leq \int_0^{+\infty}\mid \Phi_f(u) \mid \cosh(xu)du.\]
Da $a_n = \mathcal{O}(n)$ ist, haben wir
\[\mid \Phi_f(u) \mid \, \leq 2e^{2u} \cdot \sum_{n=1}^{+\infty}\mid a_n \mid e^{-2\pi n e^u} \leq 2e^{2u} \cdot \sum_{n=1}^{+\infty}n e^{-2\pi n e^u}.\]
Die Funktion
\[g(v) = v\cdot \exp(-2\pi v e^u), \quad v \geq 0, u \geq 0\]
ist monoton fallend. Nach dem Hilfssatz \ref{h6} erhalten wir die Absch\"atzung
\[2e^{2u}\cdot \sum_{n=1}^{+\infty}n\cdot \exp(-2\pi n e^u) \leq 2e^{2u} \int_0^{+\infty}v\cdot \exp(-2\pi v e^u) dv = \frac{1}{2\pi^2}e^u.\]
Daraus folgt die Behauptung.
\end{proof}
\noindent Sei $f(\tau) = \sum_{n=1}^{+\infty}a_n \exp(2\pi in\tau) \in S_2, \tau \in \mathbb{H} = \{\tau \in \mathbb{C}: \mathfrak{Im}(\tau) > 0\}$ eine Spitzenform vom Gewicht $2$. Es ist offensichtlich, dass die Asymptotik gilt
\[\Phi_f(u) = 2e^{2u}\sum_{n=1}^{+\infty}a_n e^{-2\pi n e^u} \sim 2a_1e^{2u} \cdot e^{-2\pi e^u}, \quad u \to +\infty.\]
Wir betrachten im folgenden die verf\"alschte $H_f$-Funktion
\[\hat{H}_f(z) \stackrel{defn}{=} 2a_1\int_0^{+\infty}e^{2u}e^{-2\pi e^u}\sinh(zu)du.\]
Der Begriff der verf\"alschten Funktion geht auf P\'olya zur\"uck, wo er in \cite{Polya26} die verf\"alschte $\xi$-Funktion von der Riemannschen $\xi$-Funktion betrachtet.
\begin{satz}\label{s24}
Die verf\"alschte $H_f$-Funktion $\hat{H}_f(z)$ verschwindet dort, wo $z$ die folgende Gleichung erf\"ullt
\begin{equation*}
\begin{split}
\frac{1}{(2\pi)^{z+2}}(\Gamma(z+2)-\sum_{n=0}^{+\infty}\frac{(-1)^n}{n!}\frac{(2\pi)^{z+n+2}}{z+n+2}) =  (2\pi)^{z-2}(\Gamma(2-z)-\sum_{n=0}^{+\infty}\frac{(-1)^n}{n!}\frac{(2\pi)^{-z+n+2}}{-z+n+2})
\end{split}
\end{equation*}
\end{satz}
\begin{proof}
Man setzt $t=e^u$. So hat man
\[\hat{H}_f(z) = 2a_1\int_{1}^{+\infty}te^{-2\pi t}(t^z-t^{-z})dt = I_1(z) - I_2(z)\]
mit
\[I_1(z) = 2a_1 \int_1^{+\infty}t^{z+1}e^{-2\pi t}dt, \quad I_2(z) = 2a_1 \int_1^{+\infty}t^{1-z}e^{-2\pi t}dt.\]
Man setzt $t_1 = 2 \pi t$. Dann hat man
\[I_1(z)/2a_1 = \frac{1}{(2\pi)^{z+2}}\int_{2\pi}^{+\infty}t_1^{z+1}e^{-t_1}dt_1, \quad I_2(z)/2a_1 = (2\pi)^{z-2}\int_{2\pi}^{+\infty}t_1^{1-z}e^{-t_1}dt_1.\]
Es gilt
\begin{equation*}
\begin{split}
I_1(z)/2a_1 = \frac{1}{(2\pi)^{z+2}}(\Gamma(z+2)-\int_{0}^{2\pi}t_1^{z+1}e^{-t_1}dt_1) \\ = \frac{1}{(2\pi)^{z+2}}(\Gamma(z+2)-\int_0^{2\pi}\sum_{n=0}^{+\infty}\frac{(-1)^n}{n!}t_1^{z+n+1}dt_1) = \frac{1}{(2\pi)^{z+2}}(\Gamma(z+2)-\sum_{n=0}^{+\infty}\frac{(-1)^n}{n!}\frac{(2\pi)^{z+n+2}}{z+n+2})
\end{split}
\end{equation*}
Analog
\[I_2(z)/2a_1 = (2\pi)^{z-2}(\Gamma(2-z)-\sum_{n=0}^{+\infty}\frac{(-1)^n}{n!}\frac{(2\pi)^{-z+n+2}}{-z+n+2}).\]
Damit wurde der Satz bewiesen.
\end{proof}
\begin{satz}\label{s25}
Die verf\"alschte $H_f$-Funktion $\hat{H}_f(z)$ hat nur rein-imagin\"are Nullstellen und zwar unendlich viele davon.
\end{satz}
\begin{proof}
Durch partielle Integration erhalten wir
\[\int_1^{+\infty}t^{z+1}e^{-2\pi t}dt = \frac{e^{-2\pi}}{2\pi} + \frac{z+1}{2\pi}\int_1^{+\infty}t^z e^{-2\pi t}dt\]
und
\[\int_1^{+\infty}t^{1-z}e^{-2\pi t}dt = \frac{e^{-2\pi}}{2\pi} + \frac{1-z}{2\pi}\int_1^{+\infty}t^{-z}e^{-2\pi t}dt.\]
Die Funktion $\hat{H}_f(z)$ verschwindet dort, wo $z = x + iy$ die Gleichung 
\[\sqrt{(x+1)^2 + y^2}\cdot \big \vert \int_1^{+\infty}t^z e^{-2\pi t} dt \big \vert \, = \, \sqrt{(x-1)^2+y^2} \cdot \big \vert \int_1^{+\infty}t^{-z} e^{-2\pi t} dt \big \vert\]
erf\"ullt. Die Funktion
\[F(z) = (z+1)\int_1^{+\infty}t^z e^{-2\pi t}dt - (1-z)\int_1^{+\infty}t^{-z}e^{-2\pi t}dt\]
ist eine transzendente ganze Funktion, weil das Integral $\int_1^{+\infty}t^ze^{-2\pi t}dt$ auf der ganzen komplexen Ebene $\mathbb{C}$ konvergiert. Das Integral $\int_1^{+\infty}t^ze^{-2\pi t}dt$ ist ebenfalls eine ganze Funktion. Nach dem Satz der Faktorisierung der ganzen Funktionen von Hadamard-Weierstra\ss \, hat $F(z)$ unendliche viele Nullstellen. Ist $z = x + iy$ eine Nullstelle von $F(z)$, so muss es gelten
\[\sqrt{\frac{(x+1)^2+y^2}{(x-1)^2+y^2}} = \frac{\int_1^{+\infty}t^{-x}e^{-2\pi t}dt}{\int_1^{+\infty}t^{x}e^{-2\pi t}dt}.\]
F\"ur $x>0$ ist die linke Seite $<1$ und die rechte Seite $>1$. F\"ur $x<0$ ist die linke Seite $>1$ und die rechte Seite $<1$. Daher ergibt sich $x=0$. Damit wurde der Satz bewiesen.
\end{proof}
\noindent Im folgenden wollen wir das asymptotische Verhalten von der Funktion $H_f(z)$ auf dem kritischen Streifen der Hecke-$L$-Funktion ermitteln, indem wir mit funktionentheoretischen Methoden wie in \cite[4.14]{Titch51} vorgehen. 
\begin{hilfssatz}\label{h7}
Sei $F(z)$ eine meromorphe Funktion und $\gamma$ ein Kontur, der die Nullstellen von $\sin(\pi z)$ bei $z=\rho,\rho+1,\dots,n$ enth\"alt. Falls $F(z)$ und $\sin(\pi z)$ verschiedene Polstellen besitzen, dann hat man die Formel
\[\sum_{m=\rho}^n F(m) = \frac{1}{2\pi i}\ointctrclockwise_{\gamma}\pi \cot(\pi z)F(z)dz - \sum{}^{\prime} \mathrm{Res}[\pi \cot(\pi z)F(z)],\]
wobei $\sum{}^{\prime}$ die Summe \"uber allen Polstellen von $F(z)$ innerhalb dem Kontur $\gamma$ ist.
\end{hilfssatz}
\noindent F\"ur eine Modulform $f(\tau) = \sum_{n=0}^{+\infty}a_n \exp(2\pi i n \tau), \tau \in \mathbb{H} = \{z \in \mathbb{C}: \mathfrak{Im}(z) > 0\}$ k\"onnen wir die Fourier-Koeffizienten $\{a_n\}_{n \in \mathbb{N}}$ als eine komplexe Funktion $a(z): \mathbb{C} \rightarrow \mathbb{C}$ auffassen, wobei
\begin{equation*}
a(z) = \begin{cases} a_n, \quad z = n \\ 0, \quad \text{sonst} \end{cases}
\end{equation*}
Seien $\sigma > 2$ und $f \in S_2$ eine Spitzenform vom Gewicht 2. Die Hecke-$L$-Funktion ist $L(f;s) = \sum_{n=1}^{+\infty}\frac{a_n}{n^s} = \sum_{n \leq x} \frac{a_n}{n^s} + \sum_{n > x}\frac{a_n}{n^s}, s = \sigma + it$. Der Hilfssatz \ref{h7} zeigt uns, dass es gilt
\begin{equation*}
\begin{split}
L(f;s) - \sum_{n<x}\frac{a_n}{n^s} = -\frac{1}{2i}\int_{x-i\infty}^{x+i\infty}a(z)\cot(\pi z) \frac{dz}{z^s} \\ = -\frac{1}{2i}\int_{x-i\infty}^x a(z) (\cot(\pi z) - i)\frac{dz}{z^s} - \frac{1}{2i}\int_{x}^{x+i\infty}a(z) (\cot(\pi z) + i)\frac{dz}{z^s} - a_{[x]}\frac{x^{1-s}}{1-s}
\end{split}
\end{equation*}
F\"ur $z = x + iy$ haben wir 
\[\mid \cot(\pi z) + i \mid \, = \, \frac{2}{1+e^{2 \pi y}} \leq 2 e^{-2\pi y}\]
und 
\[\mid z^{-s} \mid \, = \, z^{-\sigma}e^{t \mathrm{arg}(z)} < x^{-\sigma}e^{\mid t \mid \arctan(y/x)} < x^{-\sigma}e^{\mid t \mid \frac{y}{x}}.\]
Daher haben wir eine Absch\"atzung
\[\big \vert \int_{x}^{x+i\infty}a(z) (\cot(\pi z) + i)\frac{dz}{z^s} \big \vert \, \leq \, \mid a_{[x]} \mid x^{-\sigma}\int_0^{+\infty}e^{-2\pi y + \mid t \mid \frac{y}{x}}dy = \frac{\mid a_{[x]} \mid \cdot x^{-\sigma}}{2\pi - \frac{\mid t \mid}{x}}.\]
Analog haben wir
\[\big \vert \int_{x}^{x+i\infty}a(z) (\cot(\pi z) - i)\frac{dz}{z^s} \big \vert \, \leq \, \mid a_{[x]} \mid x^{-\sigma}\int_0^{+\infty}e^{-2\pi y + \mid t \mid \frac{y}{x}}dy = \frac{\mid a_{[x]} \mid \cdot x^{-\sigma}}{2\pi - \frac{\mid t \mid}{x}}.\]
Damit haben wir den folgenden Satz bewiesen, wobei wir nur f\"ur $f \in S_2$ formulieren, obwohl der Beweis f\"ur allgemeine Spitzenformen $f \in S_k$ auf dem gleichen Prinzip basiert. 
\begin{satz}\label{s26}
Sei $f \in S_2$ mit $L(f;s) = \sum_{n=1}^{+\infty}\frac{a_n}{n^s}, a_n = \mathcal{O}(n), s = \sigma + it$. Man hat eine asymptotische Formel gleichm\"a\ss ig auf $\sigma \geq \sigma_0 > 1, \mid t \mid \leq 2 \pi x/C$ mit einer Konstante $C>1$
\[L(f;s) = \sum_{n \leq x}\frac{a_n}{n^s}-a_{[x]}\frac{x^{1-s}}{1-s} + \mathcal{O}(x^{-\sigma+1}).\]
\end{satz}
\noindent Im folgenden schreiben wir $\log_+(x) \stackrel{defn}{=} \log(2+\mid x\mid)$. 
\begin{satz}\label{s27}
Seien $f \in S_2$ eine Spitzenform vom Gewicht $2$ und $C>1$ eine Konstante mit $\sigma > 1, \mid t \mid \leq 2 \pi x/C$. Dann hat man eine asymptotische Formel
\[H_f(z) = \mathcal{O}(\exp(-\frac{\pi}{2}\mid t \mid + \mathcal{O}((\sigma+1)\log_+(\mid\sigma+1\mid+\frac{2\pi x}{C})))\cdot \max\{\zeta(\sigma),\frac{x^{-\sigma}}{\sigma},x^{1-\sigma}\}).\]
\end{satz}
\begin{proof}
Wir haben die Stirlingsche Formel
\[\Gamma(1+\sigma+it) = \exp((\sigma+1+it-1/2)\log(\sigma+1+it)-(\sigma+1+it)+\log\sqrt{2\pi} + \mathcal{O}(\frac{1}{\mid \sigma+1+it \mid})),\]
wo man den Hauptzweig von dem komplexen Logarithmus verwendet. Insbesondere, es gilt dann
\[\Gamma(1+\sigma+it) = \mathcal{O}(\exp((\sigma+1/2)\log\mid\sigma+it\mid-t\arctan(\frac{t}{\sigma})-\sigma)).\]
Da es
\[\arctan(\frac{t}{\sigma}) = \frac{\pi}{2}\cdot \mathrm{sgn}(t)+\mathcal{O}(\frac{\sigma}{\sigma+\mid t \mid})\]
gilt, hat man dann
\[\Gamma(1+\sigma+it) = \mathcal{O}(\exp(-\frac{\pi}{2}\mid t \mid + \mathcal{O}((\sigma+1)\log_+(\mid \sigma+1\mid + \mid t \mid))).\]
Die Behauptung folgt nun aus dem Satz \ref{s26}.
\end{proof}

\end{document}